\documentclass[11pt,oneside]{amsart}
\usepackage{graphicx}
\usepackage{psfrag}
\usepackage{amscd}

\usepackage{epsfig}

\usepackage{amsmath,ifthen, amsfonts, amssymb, amsopn,amsthm}
\usepackage{color}

\catcode`\@=11
\long\def\@savemarbox#1#2{\global\setbox#1\vtop{\hsize\marginparwidth
  \@parboxrestore\tiny\raggedright #2}}

\newcommand {\me}{\medskip}
\newcommand {\sm}{\smallskip}

\newcommand{\showcomments}{yes}

\newsavebox{\commentbox}
{\ifthenelse{\equal{\showcomments}{yes}}%
{\footnotemark
        \begin{lrbox}{\commentbox}
        \begin{minipage}[t]{1.25in}\raggedright\sffamily\tiny
        \footnotemark[\arabic{footnote}]}
{\begin{lrbox}{\commentbox}}}%
{\ifthenelse{\equal{\showcomments}{yes}}%
{\end{minipage}\end{lrbox}\marginpar{\usebox{\commentbox}}}
{\end{lrbox}}}

\newcommand{\fp}{{\mathfrak p}}
\newcommand{\fs}{{\mathfrak s}}

\newcommand{\nn}{{\mathcal N}}
\newcommand{\pp}{{\mathcal P}}
\newcommand{\onn}{{\mathcal N}}
\newcommand{\bb}{{\mathcal B}}
\newcommand{\rr}{{\mathcal R}}
\newcommand{\uu}{{\mathcal U}}

\newtheorem{thm}{Theorem}[section]
\newtheorem{lem}[thm]{Lemma}

\newtheorem{cor}[thm]{Corollary}

\newtheorem{prop}[thm]{Proposition}

\theoremstyle{definition}
\newtheorem{defn}[thm]{Definition}
\newtheorem{rem}[thm]{Remark}
\newtheorem{exmp}[thm]{Example}
\newtheorem{exmps}[thm]{Examples}

\newtheorem{notation}[thm]{Notation}
\newtheorem{conv}[thm]{Convention}

\newtheorem{quest}[thm]{Question}

\newcommand{\field}[1]{\mathbb{#1}}
\newcommand{\integers}{\ensuremath{\field{Z}}}

\newcommand{\naturals}{\ensuremath{\field{N}}}
\newcommand{\reals}{\ensuremath{\field{R}}}

\newcommand{\hyperbolic}{\ensuremath{\field{H}}}
\newcommand{\pdist}{{\mathrm{pdist}}}
\newcommand{\dist}{{\mathrm{dist}}}
\newcommand{\card}{{\mathrm{card}}}
\newcommand {\N}{\mathbb{N}} %% positive integers
\newcommand {\R}{\mathbb{R}} %% reals
\newcommand {\C}{\mathbb{C}} %% complex
\newcommand {\free}{\mathbb{F}} %% free
\newcommand {\ww}{\mathcal{W}} %% W caligr.
\newcommand {\hh}{\mathcal{H}} %%H caligr.
\newcommand {\mm}{\mathcal{M}} %%M caligr.

\newcommand {\iv}{^{-1}}

\newcommand{\norm}[1]{\ensuremath{\vert \vert #1\vert  \vert}}

\newcommand {\mdk}{measure definite kernel }

\begin{document}

\title[Kazhdan and Haagerup properties from the median viewpoint]
{Kazhdan and Haagerup properties from the median viewpoint}

\author[I.~Chatterji]{Indira Chatterji}
           \address{
The Ohio-State University\\
231 W18th Ave\\
43210 Columbus OH\\
USA.}
           \email{indira@math.ohio-state.edu}

           \author[C.~Dru\c{t}u]{Cornelia Dru\c{t}u}\thanks{  {The research of the second author was supported in part by
the ANR project ``Groupe de recherche de G\'eom\'etrie et Probabilit\'es dans
les Groupes''.}}
           \address{
Mathematical Institute\\
24-29 St Giles\\
Oxford OX1 3LB\\
United Kingdom.}
           \email{drutu@maths.ox.ac.uk}

\author[F.~Haglund]{Fr\'{e}d\'{e}ric Haglund}
           \address{
Laboratoire de Math\'{e}matiques\
Universit\'{e} de Paris XI (Paris-Sud)\\
91405 Orsay\\
France.}
           \email{frederic.haglund@math.u-psud.fr}

\date{\today}

\begin{abstract}
We prove the existence of a close connection between spaces with measured walls and median metric
spaces.

We then relate properties (T) and Haagerup (a-T-me\-na\-bi\-lity) to actions on median spaces and on
spaces with measured walls. This allows us to explore the relationship between the classical properties
(T) and Haagerup and their versions using affine isometric actions on $L^p$-spaces. It also allows us
to answer an open problem on a dynamical characterization of property (T), generalizing results of
Robertson-Steger.
\end{abstract}

\maketitle

\tableofcontents

\section{Introduction.}

\subsection{Median spaces and spaces with measured walls}

 {The theory of median spaces involves a blend of geometry,  graph theory and functional analysis.} A \emph{median space} is a metric space for which, given any triple of points, there exists a unique
\emph{median point}, that is a point which is simultaneously between any two points in that triple. A
point $x$ is said to be \emph{between} two other points $a,b$ in a metric space $(X, \dist )$ if $\dist
(a,x)+\dist (x,b)=\dist (a,b)$.

Examples of median spaces are real trees, sets of vertices of simplicial trees, $\reals^n,\, n\geq 1
,\, $ with the $\ell^1$-metric, CAT(0) cube complexes with the cubes endowed with the $\ell^1$-metric,
0-skeleta of such CAT(0) cube complexes. According to Chepoi \cite{Chepoi:graphs} the class of
1-skeleta of CAT(0) cube complexes coincides with the class of median graphs (simplicial graphs whose
combinatorial distance is median). See also \cite{Gerasimov:semisplittings} and
\cite{Gerasimov:fixedpoint} for an equivalence in the same spirit.  {Median graphs are much studied in graph theory and in computer science \cite{ChepoiBandelt} and are relevant in optimization theory (see for instance \cite{MMR} and \cite{Wildstrom} and references therein for recent applications).}

 {Median metric spaces can thus be seen as
non-discrete generalizations of 0-skeleta of CAT(0) cube complexes (and geodesic median spaces can be seen as
non-discrete generalizations of
1-skeleta of CAT(0) cube complexes), same as real trees are non-discrete generalizations of simplicial
trees.}

The ternary algebra naturally associated with a median space is called a \emph{median algebra}. There
is an important literature studying median algebras. Without attempting to give an exhaustive list, we
refer the reader to \cite{Sholander1}, \cite{Sholander2}, \cite{Nieminen}, \cite{Isbell},
\cite{BandeltHedlikova}, \cite{VandeVel:book}, \cite{Basarab}. More geometrical studies of median
spaces were started in \cite{Roller:median} or \cite{NicaMaster}. In this article, we emphasize that a
median space has a richer structure than the algebraic one, and we use this structure to study groups.

Nica in \cite{Nica:cubulating}, and Niblo and the first author in \cite{ChatterjiNiblo} proved
simultaneously and independently an equivalence of category between CAT(0) cube complexes, possibly
infinite dimensional, and \emph{discrete spaces with walls} (notion introduced by F. Paulin and the
third author in \cite{HaPa}). Discrete spaces with walls were generalized by Cherix, Martin and Valette
in \cite{CherixMartinValette} to \emph{spaces with measured walls}.  {Spaces with measured walls are naturally endowed with a (pseudo-)metric.
It turns out (Corollary \ref{cor:premedian}) that a (pseudo-)metric on a space is induced by a structure of measured walls if and only if it is induced by an embedding of the space into a median space (i.e. it is submedian, in the sense of Definition \ref{defn:submed}). This is a consequence of the following results:}
\begin{thm}\label{res}
\begin{enumerate}
\item Any space with measured walls $X$ embeds isometrically in a canonically associated median space $\mathcal{M} (X)$.
Moreover, any homomorphism between two spaces with measured walls induces an isometry between the
associated median spaces.

\sm

\item Any median space $(X,\dist)$ has a canonical structure of space with measured
walls, and the wall metric coincides with the original metric. Moreover, any isometry between median
spaces induces an isomorphism between the structures of measured walls.

\sm

\item Any median space $(X,\dist)$ embeds isometrically in $L^1(\ww,\mu)$, for some measured space
$(\ww,\mu)$.
\end{enumerate}\end{thm}

 {The fact that a median space embeds into an $L^1$--space was known previously, though the embedding was not explicitly constructed, but obtained via a result of Assouad that a metric space embeds into an $L^1$--space if and only if any finite subspace of it embeds (\cite{AssouadDeza}, \cite{Assouad:L1}, \cite{Assouad:analytiqueComment}, \cite{Verheul:book}). It is moreover known that complete median normed spaces are linearly isometric to $L^1$--spaces \cite[Theorem III.4.13]{Verheul:book}.}

We note here that there is no hope of defining a median space containing a space with measured walls
and having the universality property with respect to embeddings into median spaces (see Remark
\ref{nofunctor}). Nevertheless, the medianization $\mathcal{M} (X)$ of a space with measured walls $X$
appearing in Theorem \ref{res}, (1), is canonically defined and it is minimal in some sense. This is
emphasized for instance by the fact that, under some extra assumptions, a space with measured walls $X$
is at finite Hausdorff distance from $\mathcal{M} (X)$ \cite{CDH-geometry}. In particular, it is the
case when $X$ is the $n$-dimensional real hyperbolic space with the natural structure of space with
measured walls (see Example \ref{exmp:realhyper}).

\medskip

\subsection{Properties (T) and Haagerup, actions on $L^p$-spaces.}

Topological groups with Kazhdan's property (T) act on real trees with bounded orbits, moreover with
global fixed point if the tree is complete (\cite[Corollary 5.2]{Bozejko:trees}, see also \cite{Alperin:propT}, \cite{Watatani}, \cite[Theorems 2.10.4 and 2.12.4]{BekkaHarpeValette}). The converse implication however does not hold
in general. Coxeter groups with every pair of generators satisfying a non-trivial relation act on any
real tree with fixed point (an application of Helly's Theorem in real trees); on the other hand these
groups are known to have the Haagerup property, also called a-T-me\-na\-bi\-lity \cite{BozejkoJS}.
Nevertheless, if one extends the bounded orbits property from actions on real trees to actions on
median spaces, the equivalence with property (T) does hold, for locally compact second countable
groups. More precisely, Theorem \ref{res} allows us to prove the following median characterizations of
property (T) and Haagerup property:
\begin{thm}\label{impl1}
Let $G$ be a locally compact second countable group.
\begin{enumerate}
\item\label{1-t} The group $G$ has property
(T) if and only if any continuous action by isometries on a median space has bounded orbits.
\item\label{2-a-t} The group $G$ has the Haagerup property if and only if it admits a proper continuous
 action by isometries on a median space.
\end{enumerate}
\end{thm}

 {Note that the direct implication in (\ref{1-t}) and the converse implication in (\ref{2-a-t}) follow immediately from known results of Delorme-Guichardet, respectively Akemann-Walter (see Theorem \ref{THker} in this paper), and from the fact that median spaces embed into $L^1$--spaces (see for instance \cite[Theorem V.2.4]{Verheul:book}).}  {For discrete countable groups geometric proofs of the same implications are provided implicitly in \cite{NibloReeves:groupsactingCATO}, \cite{NibloRoller:cubesandT} (see also \cite{Roller:median}) and explicitly in \cite{NicaMaster}. Nica conjectured \cite{NicaMaster} that the converse implication in (\ref{1-t}) and the direct implication in (\ref{2-a-t}) hold for discrete countable groups. This is answered in the affirmative by Theorem \ref{impl1}.}

\medskip

Theorem \ref{impl1} can be reformulated in terms of spaces with measured walls as follows:
\begin{thm}\label{implmw}
Let $G$ be a locally compact second countable group.
\begin{enumerate}\item The group $G$ has property
(T) if and only if any continuous action by automorphisms on a space with measured walls has bounded
orbits (with respect to the wall pseudo-metric).
\item The group $G$  has the Haagerup property if and only if it admits a proper continuous action
by automorphisms on a space with measured walls.
\end{enumerate}
\end{thm}

The equivalence in Theorem \ref{implmw} improves the result of Cherix, Martin and Valette
\cite{CherixMartinValette}, who showed the same equivalence for discrete groups.

Using Theorem \ref{implmw}, the classical properties (T) and Haagerup can be related to their versions
for affine actions on $L^p$-spaces.

\begin{defn}\label{atp}
Let $p>0$, and let $G$ be a topological group.
\begin{itemize}
    \item[(1)] The group $G$ \emph{has property }$FL^p$ if any affine isometric continuous action of $G$ on a space $L^p(X,\mu )$ has bounded
    orbits (equivalently, for $p>1$, it has a fixed point).

    \me

    \item[(2)] The group $G$ is \emph{a}-$FL^p$-\emph{menable} if it has a proper affine isometric continuous action on
    some space $L^p(X,\mu )$.
\end{itemize}
\end{defn}

Property $FL^2$ is equivalent to property $FH$, i.e. the fixed point property for affine continuous
actions on Hilbert spaces, and the latter is equivalent to property (T) for $\sigma$-compact groups, in
particular for second countable locally compact groups, as proved in \cite{Guichardet:bull} and
\cite{Delorme}. Likewise a-$FL^2$-menability is equivalent to a-T-menability (or Haagerup property).

Theorem \ref{implmw} and a construction from \cite{CherixMartinValette} and \cite{CornTessValette:lp}
associating to every action on a space with measured walls an affine isometric action on an $L^p$-space
implies the following.

\begin{cor}\label{c3}
Let $G$ be a second countable locally compact group.
\begin{enumerate}
    \item\label{p-t} If $G$ has property $FL^p$ for some $p>0$
    then $G$ has property (T).
    \item\label{p-a-t} If $G$ has the Haagerup property then for every $p>0$ the group $G$ is
    a-$FL^p$-menable.
\end{enumerate}
\end{cor}

Cornulier, Tessera and Valette proved the implication in (\ref{p-a-t}) for countable discrete groups
\cite[Proposition 3.1]{CornTessValette:lp}.

The implication in (\ref{p-a-t}) with $p\in (1,2)$ and a proper action on $L^p([0,1])$ has been
announced in \cite{PNowak:lp}, and a complete proof has been provided in \cite{PNowak:arxiv} (see Remark \ref{nowak}).

\begin{rem}\label{p02}
The converse statements in Corollary \ref{c3} hold for $p\in (0,2]$, in the following strengthened
version:
\begin{enumerate}
    \item\label{it}  If $G$ has property (T) then it has property $FL^p$ for every $p\in (0,2]$.

    \sm

    \item\label{iat} The group $G$ has the Haagerup property if it is a-$FL^p$-menable for some  $p\in (0,2]$.
\end{enumerate}

These statements (even slightly generalized, see Corollary \ref{c2}) follow from results of
Delorme-Guichardet (\cite{Guichardet:bull}, \cite{Delorme}) and Akemann-Walter \cite{AkermannWalter},
combined with a classical Functional Analysis result \cite[Theorem 4.10]{WellsWilliams:Embeddings}.

Thus, property (T) is equivalent to all properties $FL^p$ with $p\in (0,2]$. Likewise, the Haagerup
property is equivalent to a-$FL^p$-menability for every $p\in (0,2]$. For a discussion of the cases
when $p>2$ see Section \ref{current}.
\end{rem}

\medskip
We prove Theorems \ref{impl1} and \ref{implmw} using median definite kernels, that turn out to coincide
with the Robertson-Steger measure definite kernels \cite{RobertsonSteger:negdefker} and are very
natural settings for these notions. Along the way this allows us to answer Robertson-Steger question
\cite[Question (i)]{RobertsonSteger:negdefker} whether measure definite kernels can be given an
intrinsic characterization among the conditionally negative definite kernels (Corollary
\ref{mesdefintr}).

We also generalize to locally compact second countable groups Robertson-Steger's dynamical
characterization of property (T) \cite{RobertsonSteger:negdefker}.  This answers Open Problem 7 in
\cite{Oberwolfach}. Moreover, we give the a-T-me\-na\-bi\-lity a dynamical characterization as well.

{\begin{thm} \label{Tdyn}Let $G$ be a locally compact second countable group.
\begin{enumerate}
    \item The group $G$ has property (T) if and only if for every
measure-preserving action of $G$ on a measure space $(X , \bb , \mu )$ and every set $S \subset X$ such
that for all $g\in G$, $\mu (S\vartriangle gS) < \infty $ and $\lim_{g\to 1} \mu (S\vartriangle gS)
=0$, the supremum $\sup_{g\in G} \mu (S\vartriangle gS)$ is finite.

\me

    \item  The group $G$ is a-T-menable if and only if there exists a
measure-preserving action of $G$ on a measure space $(X , \bb , \mu )$ and there exists a set $S
\subset \bb$ such that for all $g\in G$, $\mu (S\vartriangle gS) < \infty $ and $\lim_{g\to 1} \mu
(S\vartriangle gS) =0$, but $\mu (S\vartriangle gS) \to \infty$ when $g\to \infty$.
\end{enumerate}
\end{thm}}
\subsection{Current developments and open questions}\label{current}

It is natural to ask wether the equivalence between properties (T) and $FL^p$ (respectively between
a-T-menability and a-$FL^p$-menability) can be extended to $p>2$. In \cite[$\S 3.c$]{BFGM}, by an
argument attributed to D. Fisher and G. Margulis, it is proved that for every group $G$ with property
(T) there exists $\varepsilon = \varepsilon (G)$ such that the group has property $FL^p$ for every
$p\in [1,2+\varepsilon )$. Nevertheless, no positive uniform lower bound for $\varepsilon (G)$ is
known.

For $p\gg 2$ the statements in Corollary \ref{c3} cannot be turned into equivalences. Indeed, it
follows from results of P. Pansu \cite{Pansu:cohomologie} that the group $G=Sp(n,1)$ does not have
property $FL^p$ for $p>4n+2$. More recently, Y. de Cornulier, R. Tessera and A. Valette proved in
\cite{CornTessValette:lp} that any simple algebraic group of rank one over a local field is
a-$FL^p$-menable for $p$ large enough. In particular, $G=Sp(n,1)$ is a-$FL^p$-menable for $p>4n+2$.

Also, results of M. Bourdon and H. Pajot \cite{BourdonPajot:cohomologie} imply that non-elementary
hyperbolic groups have fixed-point-free isometric actions on $\ell^p (G)$ for $p$ large enough, hence
do not have property $FL^p$. G. Yu later proved \cite{Yu:hyplp} that every discrete hyperbolic group
$G$ is a-$FL^p$-menable for $p$ large enough. In particular this holds for hyperbolic groups with
property (T).

The above quoted results of  Y. de Cornulier, R. Tessera and A. Valette, and of G. Yu, illustrate that
neither of the two converse implications in Corollary \ref{c3} hold for $p\gg 2$. This shows that for
every $p>2$ property $FL^p$ is \emph{a priori} stronger than property (T). Also, the property of
 a-$FL^p$-me\-na\-bi\-lity is a weaker version of
a-T-me\-na\-bi\-lity/Haagerup property.

\begin{quest}
Can Corollary \ref{c3} be generalized to: ``for every $p\geq q\geq 2$ property $FL^p\, $ implies
property $FL^q\, $ and a-$FL^q$-menability implies a-$FL^p$-menability'' ?
\end{quest}

\begin{quest}
Are different properties $FL^p$ and $FL^q\, $ with $p,q>2$ large enough, equivalent ? Is it on the
contrary true that for any $p_0\geq 2$ there exist groups that have property $FL^p$ for $p\leq p_0$ and
are a-$FL^p$-menable for $p>p_0$?
\end{quest}

\begin{quest}
What is the relation between $FL^p$ with $p\gg 2$ and other strong versions of property (T) defined in
terms of uniformly convex Banach spaces, like for instance the one defined in
\cite{Lafforgue:Trenforce}~?
\end{quest}

Note that, like other strong versions of property (T), the family of properties $FL^p$ separates the
semisimple Lie groups of rank one from the semisimple Lie groups with all factors of rank at least $2$
(and their respective lattices). By the results of G. Yu \cite{Yu:hyplp} all cocompact rank one
lattices are a-$FL^p$-menable for $p$ large enough. On the other hand, lattices in semisimple Lie
groups of higher rank have property $FL^p$ for all $p\geq 1$, by results of Bader, Furman, Gelander and
Monod \cite{BFGM}.

Note also that the other possible version of property (T) in terms of $L^p$-spaces, namely that
``almost invariant vectors imply invariant vectors for linear isometric actions'', behaves quite
differently with respect to the standard property (T); namely the standard property (T) is equivalent
to this $L^p$ version of it, for $1<p<\infty$ \cite[Theorem A]{BFGM}. This shows in particular that the
two definitions of property (T) (i.e. the fixed point definition and the almost invariant implies
invariant definition) are no longer equivalent in the setting of $L^p$ spaces.

\bigskip

According to Bass-Serre theory, a group splits if and only if it acts non trivially on a simplicial
tree. This implies that amalgamated products do not have property (T). Splittings were later extended
to semi-splittings, using CAT(0) cube complexes. M. Sageev showed that if $G$ is a finitely generated
group acting on a finite dimensional CAT(0) cube complex without a fixed point then there exists a
stabilizer $H$ of some convex wall such that $e(G,H)
> 1$ (see \cite{Sageev:ends}). Then, in  \cite{Gerasimov:semisplittings} and \cite{Gerasimov:fixedpoint}, V. Gerasimov removed the finite dimension assumption. (Here $e(G,H)$ stands for the number of ends of the group $G$ with
respect to the subgroup $H$, in the sense of \cite{Houghton:ends}). Conversely,  V. Gerasimov showed
that any group $G$ that has a subgroup $H$ with $e(G,H) > 1$ acts on a CAT(0) cube complex without a
fixed point, so that $H$ is a finite index subgroup in the stabilizer of a convex wall (see
\cite{Gerasimov:semisplittings}, \cite{Gerasimov:fixedpoint} and also \cite{Sageev:ends},
\cite{NibloSageevScottSw}).

Under certain stability assumptions, a non-trivial action of a group on a real tree leads to a
non-trivial action on a simplicial tree and to a splitting of the group (according to Rips,
Bestvina-Feighn \cite[Theorem 9.5]{BestvinaFeighn:stableactions}, Sela \cite[Section 3]{Sela:acces},
Guirardel \cite{Guirardel:trees}).

\begin{quest}\label{qrips}
A group $G$ acts non trivially on a median space (equivalently, $G$ does not have property (T)). Under
what assumptions is there a non trivial action of $G$ on a CAT(0) cube complex (hence a semi-splitting
of $G$ in the sense of Gerasimov-Sageev) ? On a finite dimensional CAT(0) cube complex~?
\end{quest}

Via Question \ref{qrips}, Theorem \ref{impl1} relates to one implication in M. Cowling's conjecture
(stating that a countable discrete group is a-T-menable if and only if it is weakly amenable with
Cowling-Haagerup constant 1 \cite[$\S 1.3.1$]{CCJJV}). Indeed, Guentner and Higson
\cite{GuetnerHigson:catO} showed that a countable discrete group acting properly on a finite
dimensional CAT(0) cubical complex is weakly amenable with Cowling-Haagerup constant 1. If for a
discrete countable a-T-menable group it would be possible (under extra-hypotheses) to extract from its
proper action on a median space a proper action on the 1-skeleton of a finite dimensional CAT(0)
cubical complex, then by Guentner and Higson weak ame\-na\-bi\-lity would follow. Extra-hypotheses are
needed: recent results show that the implication ``a-T-menable $\Rightarrow$ weakly amenable with
Cowling-Haagerup constant 1'' does not hold in full generality. More precisely, a wreath product $H \wr
\free_2$, where $H$ is finite and $\free_2$ is the free group on two generators, is a-T-menable
according to Cornulier-Stalder-Valette \cite{CornStaldValette:lamp}, but cannot be weakly amenable with
Cowling-Haagerup constant 1 according to Ozawa-Popa \cite[Corollary 2.11]{OzawaPopa}.

\subsection{Plan of the paper}
The paper is organized as follows. Section \ref{subs:prelmed} gives a general introduction of median
spaces (further geometric considerations on such spaces will be found in \cite{CDH-geometry}) and
proves some general results used in the sequel. In Section \ref{sec:medianization} we recall the notion
of measured wall spaces and show how those embed isometrically in a median space, proving Theorem
\ref{res} part (1). In Section \ref{sec:medianalgebras} we outline a few known results on median
algebras. We emphasize on the results needed for Section \ref{med=mws}, which explains a structure of
measured wall spaces hidden in a median space. Section  \ref{sect:kernel} is devoted to the study of several
types of kernels, leading to the proof of Theorem \ref{impl1} and its consequences.

\medskip

\noindent {\it Acknowledgments:} Part of the work on the present paper was carried out during visits to
the Universities of Paris XI (Paris-Sud) and Lille 1. The authors would like to express their gratitude
towards these institutions.
%They would also like to thank Alain Valette for useful conversations.

The first author is thankful to the FIM in Z\"urich for its hospitality, during part of the work on
this paper. The second author thanks the Centre Interfacultaire Bernoulli in Lausanne for its
hospitality during the final stages of the work on the paper. The third author wishes to thank his
colleague Alano Ancona for numerous discussions during the elaboration of this paper.

We also thank Yves de Cornulier, Pierre de la Harpe, Nicolas Monod, Guyan Robertson and Alain Valette
for useful comments.
%

%%%%%%%%%%%%%%%%%%%%%%%%%%%%%%%%%%%%%%%%%%%%%%%%%%%%%%%%%%%%%%%%%%%%%%%%%%%%%%
\section{Median spaces.}\label{subs:prelmed}
%%%%%%%%%%%%%%%%%%%%%%%%%%%%%%%%%%%%%%%%%%%%%%%%%%%%%%%%%%%%%%%%%%%%%%%%%%%%%%%

%In this section basic definitions related to the notion of median (pseudo)-metric space are provided, and some geometrical features of such a space are described.

\subsection{Definitions and examples.}

\medskip

A \emph{pseudo-metric} or \emph{pseudo-distance} on a set $X$ is a symmetric map $\pdist : X\times X
\to \R_+$ satisfying the triangle inequality. Distinct points $x\neq y$ with $\pdist (x,y)=0$ are
allowed.

A map $f:(X_1,\pdist_1)\to (X_2,\pdist_2)$  between two pseudo-metric spaces is an \emph{isometry} if
$\pdist_2(f(x),f(y))=\pdist_1(x,y)$. Note that $f$ is not necessarily injective.

A space $X$ with a pseudo-metric $\pdist$ has a canonical metric quotient $\widetilde{X}= X/\sim$
composed of the equivalence classes for the equivalence relation $x\sim y \Leftrightarrow \pdist
(x,y)=0$, endowed with the metric $\dist (\tilde{x},\tilde{y})=\pdist (x,y)$. We call $\widetilde{X}$
the \emph{metric quotient of }$X$. The natural projection map $X\to\widetilde{X}$ is an isometry.

\begin{notation}
If $x$ is a point in $X$ and $r\ge 0$ then $B(x,r)$ denotes the closed ball of radius $r$ around $x$,
that is the set $\{y\in X \; ;\; \pdist (y,x) \leq r\}$.

For every $Y\subseteq X$ and $r\ge 0$, we denote by $\nn_r(Y)$ the closed $r$-tubular neighborhood of
$Y$ in $X$, $\{ y\in X \; ;\; \pdist (y,Y)\leq r\}$.
\end{notation}

\begin{defn}[intervals and geodesic sequences]\label{defn:between}
Let $(X,\pdist)$ be a pseudo-metric space. A point $b$ is \emph{between $a$ and $c$} if $\pdist
(a,b)+\pdist (b,c)=\pdist (a,c)$.  We denote by $I(a,c)$ the set of points that are between $a$ and
$c$, and we call $I(a,c)$ the \emph{interval between $a$ and $c$}.
A \emph{path} is a finite sequence of points $(a_1,a_2,...,a_n)$. It is called a  \emph{geodesic sequence} if and only if
$$\pdist
(a_1,a_n)=\pdist (a_1,a_2)+\pdist (a_2,a_3)+\cdots +\pdist (a_{n-1},a_n)\, .$$

So $(a,b,c)$ is a geodesic sequence if and only if $b\in I(a,c)$.

\end{defn}
\begin{defn}[median point]
Let $a,b,c$ be three points of a pseudo-metric space $(X,\dist)$.  We denote the intersection $I(a,b)\cap I(b,c)\cap I(a,c)$ by $M(a,b,c)$, and we call any point in $M(a,b,c)$ a {\em
median point for $a,b,c$.} We note that $I(a,b)=\{x\in X,x\in M(a,x,b)\}$.

\end{defn}
\begin{defn}[median spaces]
A \emph{median (pseudo-)metric space} is a  (pseudo-)metric space in which for any three points $x,y,z$
the set $M(x,y,z)$ is non-empty and of diameter zero (any two median points are at pseudo-distance 0).
In particular a metric space is median if any three points $x,y,z$ have one and only one median point,
which we will denote by $m(x,y,z)$.
\end{defn}

Note that a pseudo-metric space is median if and only if its metric quotient is median.

A {\em strict median subspace} of a median pseudo-metric space $(X,\dist)$ is a subset $Y$ of $X$ such
that for any three points $x,y,z$ in $Y$, the set $M(x,y,z)$ is contained in $Y$.

A subset $Y\subset X$ is a {\em median subspace} if for any three points $x,y,z$ in $Y$, we have
$M(x,y,z)\cap Y\neq\emptyset$. Note that $Y$ is then median for the induced pseudo-metric.
An intersection of strict median subspaces is obviously a strict median subspace, thus any subset $Y\subset X$ is contained in a smallest strict median subspace, which we call
the \emph{strict median hull} of $Y$. When $X$ is a metric space median subspaces
are strict, thus we simplify the terminology to \emph{median hull}.

A \emph{homomorphism of median pseudo-metric spaces} is a map $f:X_1\to X_2$ between two median
pseudo-metric  spaces such that for any three points $x,y,z\in X_1$ we have $f(M_{X_1}(x,y,z))\subset
M_{X_2}(f(x),f(y),f(z))$. This is equivalent to asking that $f$ preserves the betweenness relation,
that is $f(I(a,b))\subset I(f(a),f(b))$.

\begin{rem}
A median metric space together with the ternary operation $(x,y,z)\mapsto m(x,y,z)$ is a particular
instance of what is called a {\em median algebra} (see Example~\ref{exmp:mmetricmalg} in
Section~\ref{sec:medianalgebras}). We will use freely some classical results in the theory of abstract
median algebras - although it is not difficult to prove them directly in our geometric context.
\end{rem}

\begin{conv}
Throughout the paper, we will call median metric spaces simply \emph{median spaces}.
\end{conv}

\begin{defn}\label{defn:submed}
We say that a metric space $(X,\dist)$  is {\em submedian} if it admits an isometric embedding into a
median space.
\end{defn}
Here are the main examples we have in mind.
\begin{exmps}\label{exmp:med}
\begin{enumerate}
\item\label{exmp:product}On the real line $\reals$, the median function is just taking the middle
point of a triple, that is $m_{\reals}(a,b,c)=a+b+c-[\max{(a,b,c)}+\min{(a,b,c)}]$.
More generally, $\reals^n$ with the $\ell_1$ norm is a median space and
$$m(\vec{x},\vec{y},\vec{z})=(m_{\reals}(x_1,y_1,z_1),\dots,m_{\reals}(x_n,y_n,z_n)).$$
The interval between two points $\vec{x},\vec{y}\in\reals^n$ is the right-angled  $n$-parallelepiped with opposite corners
$\vec{x}$ and $\vec{y}$ and edges parallel to the coordinate axes.

\item The $\ell_1$-product of two pseudo-metric spaces $(X_1,\pdist_1)$ and $(X_2,\pdist_2)$ is the set
$X_1\times X_2$, endowed with the pseudo-metric
$$\pdist((x_1,x_2),(y_1,y_2))=\pdist_1(x_1,y_1)+\pdist_2(x_2,y_2).$$
Then $(X_1\times X_2,\pdist)$ is median if and only if $(X_1,\pdist_1)$ and $(X_2,\pdist_2)$ are median
(the components of a median point in $X_1\times X_2$ are  median points of the components).
\item(trees) Every $\reals$-tree is a median space.
\item (motivating example: CAT(0) cube complexes) The $1$-skeleton of a CAT(0) cube complex is a (discrete)
median space. In fact,
    according to \cite[Theorem 6.1]{Chepoi:graphs}
    a simplicial graph is median if and only if it is the 1-skeleton of
     a CAT(0) cube complex.
\item A discrete space with walls (in the sense of \cite{HaPa}) is submedian by
\cite{ChatterjiNiblo} and \cite{Nica:cubulating}. We shall prove further in this paper that
actually a space is submedian if and only if it is a space with measured walls (see
Corollary~\ref{cor:premedian}).

\item Various examples of submedian spaces can also be deduced from Remark \ref{rem:typepsubm}. For instance, Remark \ref{rem:typepsubm} combined with results in \cite{BozejkoJS} and with Proposition
\ref{emb} implies that if $(W, S)$ is a Coxeter system and $\dist_S$ is the word distance on $W$ with
respect to $S$ then $\left( W, \dist_S^{1/2} \right)$ is submedian.

Likewise, from \cite{BallmannSwiat:polyhedra} can be deduced that if $X$ is a polygonal complex locally
finite, simply connected and of type either $(4,4)$ or $(6,3)$ and $\dist$ is its geodesic distance
then $(X, \dist^{1/2})$ is submedian.

\item($L^1$-spaces)\label{exmp:L1}  Given a measured space $(X, \mathcal{B}, \mu)$, the metric space
$L^1(X,\mu)$ is median. Indeed, it is enough to see that the real vector space ${\mathcal L}^1(X,\mu)$
of measurable functions $f:X\to\reals$ with finite $L^1$-norm is a median pseudo-metric space.
Define on ${\mathcal L}^1(X,\mu)$ a ternary operation $(f,g,h)\mapsto m(f,g,h)$ by
$$m(f,g,h)(x)=m_{\reals}(f(x),g(x),h(x)).$$
Clearly $m=m(f,g,h)$ is measurable and since it is pointwise between $f$ and $g$, it satisfies
$\norm{f-g}_1=\norm{f-m}_1+\norm{m-g}_1$. In particular $m\in {\mathcal L}^1(X,\mu)$ and $m\in I(f,g)$,
where the interval is defined with respect to the pseudo-distance $\pdist(f,g)=\norm{f-g}_1$. Similarly
we have $m\in I(g,h)$ and $m\in I(f,h)$, so that $m(f,g,h)$ is a median point for $f,g,h$.

It is easy to see that a function $p\in {\mathcal L}^1(X,\mu)$ belongs to $I(f,g)$ if and only if the
set of points $x$ such that $p(x)$ is not between $f(x)$ and $g(x)$ has measure 0. It follows that
$M(f,g,h)$ is the set of functions that are almost everywhere equal to $m(f,g,h)$, so that ${\mathcal
L}^1(X,\mu)$ is a median pseudo-metric space. We conclude that $L^1(X,\mu)$ is median because it is the
metric quotient of ${\mathcal L}^1(X,\mu)$.
\item(symmetric differences)\label{exmp:symdiff} Let $(X, \mathcal{B}, \mu)$ still denote  a measured space.
For any subset $A\subset X$, we define
$${\bb}_A=\{B\subseteq X\,|\,A\vartriangle B\in\bb\, ,\, \mu(A \vartriangle B)<+\infty\}.$$
Notice that we don't require the sets in ${\bb}_A$ to be measurable, only their symmetric difference with $A$ should be. Denote as usual by $\chi_C$   the
characteristic function of a set $C$. Then
 the map $\chi^A:{\bb}_A\to {\mathcal L}^1(X,\mu)$ defined by $B\mapsto \chi_{A \vartriangle B}$ is injective. The range of $\chi^A$ consists in the class ${\mathcal S}^1(X,\mu)$ of all
  characteristic functions of measurable subsets with finite measure. Indeed the
   preimage of $\chi_{B'}$ (with $B'\in\bb,\mu(B')<+\infty$) is the subset $B:= A \vartriangle B'$.
   Observe that the $L^1$-pseudo-distance
 between two functions $\chi_{B'}$ and $\chi_{C'}$ in ${\mathcal S}^1(X,\mu)$ is equal to $\mu(B'\vartriangle C')$. Since we have
$$(A\vartriangle B)\vartriangle(A\vartriangle C)=B\vartriangle C,$$
it follows that for any two elements $B_1,B_2\in{\bb}_A$ the symmetric difference $B_1\vartriangle B_2$ is measurable with finite
 measure, and  the pull-back of the $L^1$-pseudo-distance by the bijection
 ${\bb}_A\to {\mathcal S}^1(X,\mu)$ is the pseudo-metric $\pdist_\mu$ defined by
 $\pdist_\mu(B_1,B_2)=\mu(B_1\vartriangle B_2)$.

We claim that $({\bb}_A,\pdist_\mu)$ is a median pseudo-metric space,
or equivalently that ${\mathcal S}^1(X,\mu)$ is a median subspace of ${\mathcal L}^1(X,\mu)$.
This follows easily from the explicit formula:
$$m(\chi_A,\chi_B,\chi_C) = \chi_{(A\cup B)\cap(A\cup C)\cap(B\cup C)}=\chi_{(A\cap B)\cup(A\cap C)\cup(B\cap C)}\, .$$

Note that $I(\chi_A,\chi_B)\cap {\mathcal S}^1(X,\mu)$ is composed of the characteristic functions
$\chi_C$ such that there exists $C'\in\bb$ satisfying $ \mu(C'\vartriangle C)=0$ and $ A\cap B\subset
C'\subset A\cup B$.
 Later we will prove that any median space embeds isometrically as a median subspace of some space
  ${\mathcal S}^1(X,\mu)$ (compare with the similar result in the context of median algebras
 appearing in Corollary \ref{sigmaemb}).
%
%\item The plane $\reals^2$ with the euclidean distance is submedian as the map
%
%\begin{eqnarray*}\reals^2&\to & L^1([0,2\pi])\\
%\left(\begin{array}{c}x\\y\end{array}\right)&\mapsto& x\sin(t)+y\cos(t)
%\end{eqnarray*}
%
%is easily seen to be an isometry.
\end{enumerate}\end{exmps}
\begin{rem}\label{premedl1}
{In view of Lemma~\ref{lem:embedL1} and of Example~\ref{exmp:med}, (\ref{exmp:L1}), a metric space $(X,
\dist )$ is submedian if and only if it embeds isometrically in a space $L^1(\ww , \mu)$, for some
measured space $(\ww , \mu)$. Thus, the notion of submedian space coincides with the notion of metric
space of type $1$ as defined in \cite[Troisi\`eme partie, $\S 2$]{BretDacunhaCastelKriv}. Similarly,
submedian metric is the same thing as metric of type $1$.}
\end{rem}
\begin{rem}\label{nofunctor} (1) It is not possible in general to define for every submedian space $Y$
a median completion, that is a median space
    containing an isometric copy of $Y$, and such that any isometric embedding of $Y$ into a median
    space extends to it. This can be seen in the following example.

Let $E=\reals^7$ endowed with the $\ell_1$ norm, and let $\{ e_i\; ;\; i=1,2,...,7 \}$ be the canonical
basis. Let $Y_x$ be the set composed of the four points $A,B,C,D$ in $E$ defined by
$A=\frac{x}{2}(e_1+e_2+e_3)+(1-x)e_4$, $B=\frac{x}{2}(-e_1-e_2+e_3)+(1-x)e_5$,
$C=\frac{x}{2}(e_1-e_2-e_3)+(1-x)e_6$, $D=\frac{x}{2}(-e_1+e_2-e_3)+(1-x)e_7$, where $x\in[0,1]$.

Any two distinct points in $Y_x$ are at $\ell_1$-distance $2$. Thus all $Y_x$ with the
$\ell_1$-distance  are pairwise isometric. The median hull of $Y_x$ is composed of $Y_x$ itself and of
the eight vertices of a cube of edge length $x$ defined by $\frac{x}{2}(\pm e_1\pm e_2\pm e_3)$. Thus,
for two distinct values $x\neq x'$ the median hulls of $Y_x$ and of $Y_{x'}$ are not isometric.

Note that the median hull of $Y_0$ is the simplicial tree with five vertices, four of which are
endpoints. The median hull of $Y_1$ is the set of eight vertices of the unit cube. Consequently, it
cannot even be guaranteed that two median hulls of two isometric submedian spaces are isomorphic as
median algebras.

\medskip

(2) Given a subspace $Y$ of a median space $X$, it is in general not possible to extend an isometry of
    $Y$ to an isometry (or at least an isomorphisms of median algebras) of the median hull of $Y$.

With the same notations as in (1), the isometry $Y_0\times Y_1\to Y_0\times Y_1$ switching the points
of $Y_0$ with the points of $Y_1$ cannot be extended to the median hull of $Y_0\times Y_1$ in $E\times
E$.
\end{rem}
%%%%%%%%%%%%%%%%%%%%%%%%%%%%%%%%%%%%%%%%%%%%%%%%%%%%%%%%
\subsection{Convexity and gate property in median spaces.}
%%%%%%%%%%%%%%%%%%%%%%%%%%%%%%%%%%%%%%%%
\begin{defn}\label{conv} Let $(X,\pdist)$ denote some pseudo-metric space.
A subset $Y\subset X$ is said to be \emph{convex} if for any $a,b\in Y$ the set $I(a,b)$ is contained
in $Y$. It is \emph{quasi-convex}
 if for any
$a,b\in Y$ the set $I(a,b)$ is contained in $\onn_M(Y)$ for some $M$ uniform in $a,b\in Y$.  The
\emph{convex hull} of a subset $Y\subset X$ is the intersection of all convex subsets containing~$Y$.
\end{defn}
Note that any convex subspace of a median space is median but not the converse, as for instance any
subset of cardinality two is a median subspace, while it might not be convex. The median hull of a
subset  is contained in the convex hull, and as the example above shows the inclusion may be strict.

We now introduce  a notion which is  related to convexity in median spaces, and which is commonly used
in the theory of Tits buildings (see for example \cite{Scharlau85})  {and in graph theory (\cite{Mulder:book}, \cite{VandeVel:book}).}
\begin{defn}[gate]\label{defn:gate}
Let $(X,\dist)$ be a metric space, let $Y$ be a subset of $X$, and $x$ some point in $X$.

We say that a point $p\in X$ is {\em between $x$ and $Y$} if it is between $x$ and any $y\in Y$. When a
point $p\in Y$ is between $x$ and $Y$, we say that $p$ is a {\em gate between $x$ and $Y$}. Note that
there is always at most one gate $p$ between $x$ and $Y$, and that $\dist(x,p)=\dist(x,Y)$.

We say that $Y$ is {\em gate-convex} if for any point $x\in X$ there exists a gate (in $Y$) between $x$
and $Y$. We then denote by $\pi_Y(x)$ this gate, and call the map $\pi_Y$ the {\em projection map onto
$Y$}.
\end{defn}
%
%%%%%%%%%%%%%%%
%%%%%%%%%%%%%%%%
%%%%%%%%%%%%%%%
\begin{lem}[gate-convex subsets]\label{lem:gate}

\begin{enumerate}
\item The projection map onto a gate-convex subset is 1-Lipschitz.

\sm

\item Any gate-convex subset is closed and convex.

\sm

\item In a complete median space, any closed convex subset is gate-convex.
\end{enumerate}
\end{lem}

In other words, for closed subsets of a complete median space, convexity is equivalent to
gate-convexity.
\begin{proof} (1)
Let $x,x'$ be two points in a metric space $X$, and let $p,p'$ be the respective gates between $x,x'$
and a gate-convex subset $Y$. Since $(x,p,p')$ and $(x',p',p)$ are geodesic sequences, we have that
\begin{eqnarray*}\dist(x,p)+\dist(p,p')&\leq&\dist(x,x')+\dist(x',p')\\
\dist(x',p')+\dist(p',p)&\leq&\dist(x',x)+\dist(x,p)
\end{eqnarray*}
By summing up the two inequalities, we conclude that $\dist (p , p')\le \dist (x , x')$.

\medskip

(2) Assume that $Y$ is gate-convex and that $(x,y,z)$ is a geodesic sequence with $x,z\in Y$. Let $p$
be the gate between $y$ and $Y$, so that $(y,p,x)$ and $(y,p,z)$ are geodesic sequences. Hence
$(x,p,y,p,z)$ is a geodesic sequence, which forces $y=p\in Y$.

Any point $x$ in the closure of $Y$ satisfies $\dist(x,Y)=0$. Thus if $p$ is the gate between $x$ and
$Y$ we have $\dist(x,p)=0$, hence $x\in Y$. We conclude that $Y$ is closed.

\medskip

(3) Let $Y$ be a closed convex subset of a complete median space $X$. For any $x\in X$ choose a
sequence $(y_k)_{k\ge 0}$ of points in $Y$ such that $\dist(y_k,x)$ tends to  $\dist(x,Y)$. First
observe that $(y_k)_{k\ge 0}$ is a Cauchy sequence.
 Indeed, denote by $\epsilon_k = \dist (y_k, x)- \dist (Y,x)$, which clearly is a sequence of positive
 numbers converging to zero.
  Let $m_{k,\ell}$ be the median point  of $(x,y_k,y_\ell)$. Then
 $\dist(x,y_k)+\dist(x,y_\ell)=2\dist(x,m_{k,\ell})+\dist(y_k,y_\ell)$
 and so by convexity of $Y$ we have $\dist(x,y_k)+\dist(x,y_\ell)\ge 2\dist(x,Y)+\dist(y_k,y_\ell)$.
 It follows that $\dist(y_k,y_\ell)\leq \epsilon_k + \epsilon_\ell$.
 Since $X$ is complete the sequence $(y_k)_{k\ge 0}$ has a limit $p$ in $X$. Since $Y$ is closed,
 the point $p$ is in $Y$. Note that $\dist(x,p)=\dist(x,Y)$. It remains to check that $p$ is between $x$
and $Y$.

Let $y$ be some point in $ Y$, and let $m$ be the median point of $x,p,y$. By convexity of $Y$ we have
$m\in Y$, so that $\dist(x,m)\ge \dist(x,Y)$. We also have $\dist(x,p)=\dist(x,m)+\dist(m,p)$. Since
$\dist(x,p)=\dist(x,Y)$ we get $\dist(m,p)=0$ as desired.
\end{proof}
We now prove that in a median space the metric intervals are gate-convex.
\begin{lem}\label{lem:medprojection}
In a median metric space any interval $I(a,b)$ is gate-convex, and the gate between an arbitrary point
$x$ and $I(a,b)$ is $m(x,a,b)$.
\end{lem}
\begin{proof}
Consider an arbitrary point $x$ in the ambient median metric space $X$, $p$ the median point $m(x,a,b)$
and $y$ an arbitrary point in $I(a,b)$. We will show that $(x,p,y)$ is a geodesic sequence.

We consider the median points  $a'=m(x,a,y)$, $b'=m(x,b,y)$ and $p'=m(x,a',b')$. Note that $p'\in
I(x,a')\subset I(x,a)$ and similarly $p'\in I(x,b)$.

Since $(a,y,b)$, $(a,a',y)$ and $(y,b',b)$ are geodesic sequences, the sequence $(a,a',y,b',b)$ is
geodesic as well. So $I(a',b')\subset I(a,b)$, hence $p'\in I(a,b)$.

We proved that $p'\in I(x,a) \cap I(x,b)\cap I(a,b)$, which by the uniqueness of the median point
implies $p'=p$. It follows that $p\in I(x,a')\subset I(x,y)$.
%%% doesn't use rectangles at all !!!
\end{proof}
We can now deduce that the median map is 1-Lipschitz, in each variable and on $X\times X\times X$
endowed with the $\ell_1$-metric.
\begin{cor}\label{cor:medlip}
Let $X$ be a median space.
\begin{enumerate}
    \item For any two fixed points $a,b\in X$ the interval $I(a,b)$ is closed and convex, and the map $x\mapsto
m(x,a,b)$ is 1-Lipschitz.

\sm

    \item  The median map $m:X\times X\times X\to X$ is 1-Lipschitz (here $X\times
X\times X$ is endowed with the $\ell_1$-product metric as defined in Example~\ref{exmp:med},
(\ref{exmp:product})).
\end{enumerate}
\end{cor}
\begin{proof}
Combine Lemma~\ref{lem:medprojection} and Lemma~\ref{lem:gate}, and use the fact that, given six points
$a,b,c,a',b',c'\in X$, the distance between the median points $m(a,b,c)$ and $m(a',b',c')$ is at most
$$\dist(m(a,b,c),m(a',b,c))+\dist(m(a',b,c),m(a',b',c))+\dist(m(a',b',c),m(a',b',c'))\, .$$ \end{proof}
\subsection{Approximate geodesics and medians; completions of median spaces.}

\medskip

We prove that the median property is preserved under metric completion. In order to do it, we need an
intermediate result stating that in a median space, approximate geodesics are close to geodesics, and
approximate medians are close to medians. We begin by defining approximate geodesics and medians.

\begin{defn}\label{defn:quasigeodesic}
Let $(X,\dist)$ be a metric space and let $\delta$ be a non-negative real number. We say that {\em $z$
is between $x$ and $y$ up to $\delta$} provided
$$\dist (x,z)+\dist (z, y )\le \dist
(x,y)+\delta \, .
$$

We say that $(a_1,a_2,...,a_n)$ is a $\delta$-\emph{geodesic sequence} if
$$\dist (a_1,a_2)+\dist
(a_2,a_3)+\cdots +\dist (a_{n-1},a_n)\leq \dist (a_1,a_n) +\delta \, .
$$
\end{defn}

\begin{notation}\label{d-geod-med}
Let $x,y$ be two points of $X$. We denote by $I_\delta(a,b)$ the set of points that are between $a$ and
$b$ up to $\delta$.

Let $x,y,z$ be three points of $X$. We denote by $M_\delta(a,b,c)$ the intersection
$$I_{2\delta}(a,b)\cap I_{2\delta}(b,c)\cap I_{2\delta}(a,c)\, .$$
\end{notation}

In accordance with the previous notation, whenever $\delta=0$  the index is dropped.

\smallskip

\begin{lem}\label{intermediate}
Given $\delta , \delta'\geq 0\, ,$ for every $c\in I_\delta (a,b)$ the set $I_{\delta'}(a,c)$ is
contained in $ I_{\delta + \delta' }(a,b)$.
\end{lem}
\begin{defn}\label{defn:deltamed}
Let $x,y,z$ be three points in a metric space. If $M_\delta(x,y,z)$ is non-empty then any point in it
is called a {\em $\delta$-median point for $x,y,z$}.
\end{defn}
\begin{lem}\label{lem:medianisquasimed}
Let  $(X,\dist )$ be a median space, and $a,b,c$ three arbitrary points in it.
\begin{enumerate}
    \item[(i)] The set $I_{2\delta}(a,b)$ coincides with $\onn_\delta \left( I(a,b) \right)$.

    \me

    \item[(ii)] The following sequence of inclusions holds:
    \begin{equation}\label{eq:mdelta}
   B(m(a,b,c), \delta) \subseteq M_\delta(a,b,c)\subseteq B(m(a,b,c), 3\delta)\, .
\end{equation}
\end{enumerate}
\end{lem}
\begin{proof} Statement (i) immediately follows from Lemma \ref{lem:medprojection}.

The first inclusion in (\ref{eq:mdelta}) is obvious. We prove the second inclusion. To do it, we
consider the median points $p_1=m(p,a,b),p_2=m(p,b,c),p_3=m(p,a,c),q=m(p_1,b,c),r=m(q,a,c)$.

First we show that $r=m(a,b,c)$. Indeed $r\in I(a,c)$ by definition. We also have $r\in I(q,c)$, and
since $q\in I(c,b)$ it follows that $r\in I(b,c)$. Finally we have $r\in I(a,q)$. Now $q\in I(p_1,b)$
and $p_1\in I(a,b)$, so $q\in I(a,b)$. It follows that $r\in I(a,b)$.

It remains to estimate the distance between $p$ and $r$. According to (i) and Lemma
\ref{lem:medprojection} the point $p$ is at distance at most $\delta$ from $p_1,p_2$ and $p_3$
respectively.

By Corollary~\ref{cor:medlip} we have $\dist (p_2,q)\le \dist (p,p_1)\le\delta$. Hence $\dist (p,q)\le
2\delta$. Applying Corollary~\ref{cor:medlip} again we get $\dist (p_3,r)\le \dist (p,q)\le 2\delta$,
consequently $\dist (p,r)\le 3\delta$.
\end{proof}

 {The following result is also proved in \cite[Corollary II.3.5]{Verheul:book}. For completeness we give another proof here.}

\begin{prop}\label{cor:medcomplete}
The metric completion of a median space is a median space as well.
\end{prop}

\begin{proof}  Let $(X,\dist)$ be a median space, and let $X\to \hat X$ be
 the metric completion. For simplicity we denote the distance on $\widehat X$ also by ${\dist}$.

The median map $m:X\times X\times X\to X\subset \widehat X$ is 1-Lipschitz by
Corollary~\ref{cor:medlip}. Thus it extends to a 1-Lipschitz map $\widehat{X}\times \widehat{X}\times
\widehat{X}\to \widehat X$, also denoted by $m$.

Clearly for any three points  ${a},{b},{c}$ in $\widehat X$, the point ${m}({a},{b},{c})$ is  median
for ${a},{b},{c}$. We now prove that $m(a,b,c)$ is the unique median point for $a,b,c$.  Let $p$ be
another median point for $ a, b,c$. The points  ${a},{b},{c}$ are limits of sequences
$(a_n),(b_n),(c_n)$ of points in $X$.   Let $m_n$ be the median point of $a_n,b_n,c_n$. Set
$\delta_n=\dist (a,a_n)+\dist (b,b_n)+\dist (c,c_n)$.

We show that $p$ is a $\delta_n$-median point for $a_n,b_n,c_n$. Indeed we have that $\dist
(a_n,p)+\dist (p,b_n)$ is at most $\dist (a_n,a)+\dist (a,p)+\dist (p,b)+\dist (b,b_n) =\dist
(a_n,a)+\dist (a,b)+\dist (b,b_n)\le 2\dist (a,a_n)+\dist (a_n,b_n)+2\dist (b,b_n)\le \dist
(a_n,b_n)+2\delta_n$. The other inequalities are proved similarly.

The point $p$ is also the limit of a sequence of points $p_n$ in $X$, such that $\dist(p,p_n)\le
\delta_n$. It follows that $p_n$ is a $2\delta_n$-median point for $a_n,b_n,c_n$. By
Lemma~\ref{lem:medianisquasimed} we then have that $\dist (p_n,m_n)\le 6\delta_n$. Since $\delta_n\to
0$ we get $p=m(a,b,c)$.

\end{proof}

\subsection{Rectangles and parallel pairs.}

In a median space $X$, the following notion of rectangle will allow us to treat median spaces as a
continuous version of the 1-skeleton of a CAT(0) cube complex.

\begin{defn}\label{defn:rectangle}
A \emph{quadrilateral}  in a metric space $(X,\dist)$ is a closed path $(a,b,c,d,a)$, which we rather denote by $[a,b,c,d]$. A quadrilateral $[a,b,c,d]$ is a \emph{rectangle} if  the four sequences
$(a,b,c)$, $(b,c,d)$, $(c,d,a)$ and $(d,a,b)$ are geodesic.\end{defn}
\begin{rem}\label{rem:rect}
\begin{enumerate}
\item By the triangular inequality, in a rectangle $[a,b,c,d]$ the following equalities hold: $\dist
(a,b)=\dist (c,d)$, $\dist (a,d)=\dist (b,c)$ and $\dist
(a,c)=\dist (b,d)$.

\item(rectangles in intervals) If $x,y\in I(a,b)$ then $[x,m(x,y,a),y,m(x,y,b)]$ is a rectangle.

\item(subdivision of rectangles) Let $[a,b,c,d]$ be a rectangle.
Let $e\in I(a,d)$ and $f=m(e,b,c)$. Then $[a,b,f,e]$ and $[c,d,e,f]$  are rectangles.\end{enumerate}
\end{rem}

\begin{defn}\label{defn:parallel}(parallelism on pairs)
Two pairs $(a,b)$ and $(d,c)$ are \emph{parallel} if $[a,b,c,d]$
is a rectangle.
\end{defn}
We mention without proof the following remarkable fact that confirms the analogy with CAT(0) cube complexes:
\begin{prop}
In a median space the parallelism on pairs is an equivalence relation.
\end{prop}
We now explain how to any 4-tuple of points one can associate a rectangle.

\begin{lem}\label{RectCan}
Let $[x,a,y,b]$ be any quadrilateral  in a median space. Then there exists a unique rectangle
$[x',a',y',b']$ satisfying the following properties:
\begin{enumerate}
    \item\label{lateral} the following sequences are geodesic: $$(x,x',a',a),\, (a,a',y',y),\, (y,y',b',b),\,
    (b,b',x',x)\, ;$$
    \item\label{xy} $(a,a',b',b)$ is a geodesic sequence;
    \item\label{ab}  $(x,x',y')$ and $(y,y',x')$ are geodesic sequences.
\end{enumerate}
\end{lem}

\begin{proof}

\medskip

\noindent\emph{Existence.} Let $x'=m(x,a,b)$ and $y'=m(y,a,b)$, and let $a'=m(a,x',y')$ and $b'=m(b,x',y')$ (see
Figure \ref{fig1}). Then $[x',a',y',b']$ is a rectangle by Remark \ref{rem:rect}, (3). Properties
(\ref{lateral}) and (\ref{xy}) follow immediately from the construction, property (\ref{ab}) follows
from Lemma \ref{lem:medprojection} applied to $x$ and $y'\in I(a,b)$, respectively to $y$ and $x'\in
I(a,b)$.

\medskip

\noindent\emph{Uniqueness.} Let $[x',a',y',b']$ be a rectangle satisfying the three required properties.
Properties (\ref{lateral}), (\ref{xy}) and the fact that $[x',a',y',b']$ is a rectangle imply that
$x'=m(x,a,b)$ and $y'=m(y,a,b)$. Again property (\ref{xy}) and the fact that $[x',a',y',b']$  is a
rectangle imply that $a'=m(a,x',y')$ and $b'=m(b,x',y')$.\end{proof}

%TeXCAD Options
%\grade{\on}
%\emlines{\off}
%\epic{\off}
%\beziermacro{\on}
%\reduce{\on}
%\snapping{\off}
%\quality{8.00}
%\graddiff{0.01}
%\snapasp{1}
%\zoom{4.0000}
\unitlength 0.75mm % = 2.85pt
\linethickness{0.4pt}
\ifx\plotpoint\undefined\newsavebox{\plotpoint}\fi % GNUPLOT compatibility
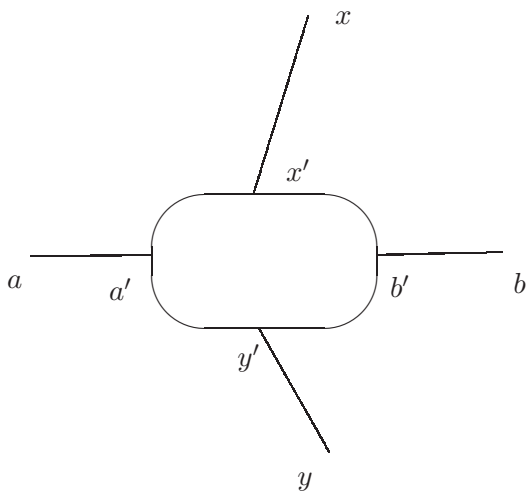
\begin{picture}(90,90)(0,0)
\put(91.75,39.13){\oval(40,23.75)[]}
%\emline(10.25,55)(31.75,55.25)
\multiput(50.25,40)(2.6875,.03125){8}{\line(1,0){2.6875}}
%\end
%\emline(72,55.25)(94,55.75)
\multiput(112,40.25)(1.466667,.033333){15}{\line(1,0){1.466667}}
%\end
%\emline(59.5,97.75)(49.75,66)
\multiput(99.5,82.75)(-.033737024,-.109861592){289}{\line(0,-1){.109861592}}
%\end
%\emline(50.75,42.25)(63.25,20.25)
\multiput(90.75,27.25)(.033692722,-.059299191){371}{\line(0,-1){.059299191}}
%\end
\put(47.5,35.5){\makebox(0,0)[cc]{$a$}} \put(66.25,34.25){\makebox(0,0)[cc]{$a'$}}
\put(115.75,34.75){\makebox(0,0)[cc]{$b'$}} \put(137,35.25){\makebox(0,0)[cc]{$b$}}
\put(105.75,82.25){\makebox(0,0)[cc]{$x$}} \put(97.75,55.25){\makebox(0,0)[cc]{$x'$}}
\put(89,22){\makebox(0,0)[cc]{$y'$}} \put(99,0){\makebox(0,0)[cc]{$y$}}
\end{picture}
\begin{figure}[!ht]
\centering \caption{Central rectangle.} \label{fig1}
\end{figure}

\begin{defn}\label{defn:centrect}
We call the rectangle $[x',a',y',b']$ described in Lemma \ref{RectCan} the \emph{central rectangle} associated with
the quadrilateral  $[x,a,y,b]$.
\end{defn}

\begin{figure}[!ht]
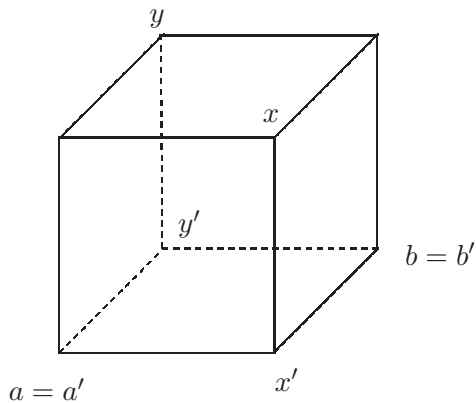

\centering \caption{Example of central rectangle.} \label{fig2}
\end{figure}

\begin{rem}\label{rem:xyab}
Property (\ref{ab}) cannot be improved to ``$(x,x',y',y)$ is a geodesic sequence'', as shown by the
example of a unit cube in $\R^3\, ,$ with $a,b$ two opposite vertices of the lower horizontal face, and
$x,y$ the two opposite vertices of the upper horizontal face that are not above $b$ or $d$ (see Figure
\ref{fig2}).

Note also that in general the central rectangle associated with $[x,a,y,b]$ is distinct from the central rectangle associated with $[a,x,b,y]$ (again see Figure~\ref{fig2}).
\end{rem}

%TeXCAD Options
%\grade{\on}
%\emlines{\off}
%\epic{\off}
%\beziermacro{\on}
%\reduce{\on}
%\snapping{\off}
%\quality{8.00}
%\graddiff{0.01}
%\snapasp{1}
%\zoom{4.0000}
\unitlength 0.75mm % = 2.85pt
\linethickness{0.4pt}
\ifx\plotpoint\undefined\newsavebox{\plotpoint}\fi % GNUPLOT compatibility
\begin{picture}(70,90)(0,0)
\put(67.5,18.5){\framebox(38,38)[cc]{}}
%\emline(45.75,136.75)(64,155)
\multiput(105.75,56.75)(.033733826,.033733826){541}{\line(0,1){.033733826}}
%\end
%\emline(7.75,136.75)(26,155)
\multiput(67.75,56.75)(.033733826,.033733826){541}{\line(0,1){.033733826}}
%\end
%\emline(45.75,98.75)(64,117)
\multiput(105.75,18.75)(.033733826,.033733826){541}{\line(0,1){.033733826}}
%\end
%\dashline{1}(7.5,98.5)(25.75,116.75)
\multiput(67.43,18.43)(.032589,.032589){20}{\line(1,0){.032589}}
\multiput(68.73,19.73)(.032589,.032589){20}{\line(1,0){.032589}}
\multiput(70.04,21.04)(.032589,.032589){20}{\line(1,0){.032589}}
\multiput(71.34,22.34)(.032589,.032589){20}{\line(1,0){.032589}}
\multiput(72.64,23.64)(.032589,.032589){20}{\line(0,1){.032589}}
\multiput(73.95,24.95)(.032589,.032589){20}{\line(0,1){.032589}}
\multiput(75.25,26.25)(.032589,.032589){20}{\line(0,1){.032589}}
\multiput(76.55,27.55)(.032589,.032589){20}{\line(0,1){.032589}}
\multiput(77.86,28.86)(.032589,.032589){20}{\line(1,0){.032589}}
\multiput(79.16,30.16)(.032589,.032589){20}{\line(1,0){.032589}}
\multiput(80.47,31.47)(.032589,.032589){20}{\line(1,0){.032589}}
\multiput(81.77,32.77)(.032589,.032589){20}{\line(1,0){.032589}}
\multiput(83.07,34.07)(.032589,.032589){20}{\line(0,1){.032589}}
\multiput(84.38,35.38)(.032589,.032589){20}{\line(1,0){.032589}}
%\end
%\dashline{1}(25.5,155)(25.75,117)
\put(85.43,74.93){\line(0,-1){.974}} \put(85.44,72.98){\line(0,-1){.974}}
\put(85.46,71.03){\line(0,-1){.974}} \put(85.47,69.08){\line(0,-1){.974}}
\put(85.48,67.13){\line(0,-1){.974}} \put(85.49,65.19){\line(0,-1){.974}}
\put(85.51,63.24){\line(0,-1){.974}} \put(85.52,61.29){\line(0,-1){.974}}
\put(85.53,59.34){\line(0,-1){.974}} \put(85.55,57.39){\line(0,-1){.974}}
\put(85.56,55.44){\line(0,-1){.974}} \put(85.57,53.49){\line(0,-1){.974}}
\put(85.58,51.55){\line(0,-1){.974}} \put(85.6,49.6){\line(0,-1){.974}}
\put(85.61,47.65){\line(0,-1){.974}} \put(85.62,45.7){\line(0,-1){.974}}
\put(85.63,43.75){\line(0,-1){.974}} \put(85.65,41.8){\line(0,-1){.974}}
\put(85.66,39.85){\line(0,-1){.974}} \put(85.67,37.9){\line(0,-1){.974}}
%\end
%\dashline{1}(25.75,117)(64,117)
\put(85.68,36.93){\line(1,0){.981}} \put(87.64,36.93){\line(1,0){.981}}
\put(89.6,36.93){\line(1,0){.981}} \put(91.56,36.93){\line(1,0){.981}}
\put(93.53,36.93){\line(1,0){.981}} \put(95.49,36.93){\line(1,0){.981}}
\put(97.45,36.93){\line(1,0){.981}} \put(99.41,36.93){\line(1,0){.981}}
\put(101.37,36.93){\line(1,0){.981}} \put(103.33,36.93){\line(1,0){.981}}
\put(105.3,36.93){\line(1,0){.981}} \put(107.26,36.93){\line(1,0){.981}}
\put(109.22,36.93){\line(1,0){.981}} \put(111.18,36.93){\line(1,0){.981}}
\put(113.14,36.93){\line(1,0){.981}} \put(115.1,36.93){\line(1,0){.981}}
\put(117.06,36.93){\line(1,0){.981}} \put(119.03,36.93){\line(1,0){.981}}
\put(120.99,36.93){\line(1,0){.981}} \put(122.95,36.93){\line(1,0){.981}}
%\end
%\emline(25.25,154.75)(63.75,154.75)
\put(85.25,74.75){\line(1,0){38.5}}
%\end
%\emline(63.75,154.75)(63.75,154.25)
\put(123.75,74.75){\line(0,-1){.5}}
%\end
%\emline(63.75,154.5)(63.75,116.75)
\put(123.75,74.5){\line(0,-1){37.75}}
%\end
\put(65.25,12){\makebox(0,0)[cc]{$a=a'$}} \put(107.75,13.75){\makebox(0,0)[cc]{$x'$}}
\put(135,36){\makebox(0,0)[cc]{$b=b'$}} \put(90.5,41.25){\makebox(0,0)[cc]{$y'$}}
\put(84.75,77.75){\makebox(0,0)[cc]{$y$}} \put(105,60.5){\makebox(0,0)[cc]{$x$}}
\end{picture}

Property (\ref{ab}) in Lemma \ref{RectCan} can be slightly improved as follows.

\begin{lem}\label{RectCan2}
Let $x,y,p,q$ be four points such that $(x,p,q)$ and $(p,q,y)$ are geodesic sequences.
 Then there exists a geodesic sequence $(x,x',y',y)$ such that $(x',y')$ and $(p,q)$ are
 parallel.
\end{lem}

\proof Applying Lemma \ref{RectCan} to the quadrilateral $[p,q,y,x]$, we note that the resulting central rectangle $[p',q',y',x']$ satisfies $p'=p,q'=q$.

\endproof

%%%%%%%%%%%%%%%%%%%%%%%%%%%%%%%%%%%%%%%%%%%%%%%%%%%%%%%%%%%%%%%%%%%%%%%%%
\section{Space with measured walls, median space associated to it.}\label{sec:medianization}
%%%%%%%%%%%%%%%%%%%%%%%%%%%%%%%%%%%%%%%%%%%%%%%%%%%%%%%%%%%%%%%%%%%%%%%

%%%%%%%%%%%%%%%%%%%%%%%%
\subsection{Preliminaries on measures.}
%%%%%%%%%%%%%%%%%%%%%%%%
We recall the relevant definitions on measured spaces. A reference is \cite{HBauer:measure}, whose
terminology we adopt here.
Let $Y$ be a non-empty set and let ${\mathcal P}(Y)$ be the power set of $Y$. A {\em ring} is a subset
of ${\mathcal P}(Y)$ containing the empty set, closed with respect to finite unions and differences.
%
%\begin{enumerate}
%\item $\emptyset\in \mathcal{R}$.
%\item If $A,B\in \mathcal{R}$ then $A\cup B\in \mathcal{R}$.
%\item If $A,B\in \mathcal{R}$ then $A\setminus B \in \mathcal{R}$.
%\end{enumerate}
%
A \emph{$\sigma$-algebra} is a subset in $\pp (Y)$ containing the empty set, closed with respect to
countable unions and containing together with any subset its complementary.
%\begin{enumerate}
%   \item $\emptyset\in \mathcal{A}$.
%   \item If $A_n\in \mathcal{B}\, ,\, n\in \N$, then $\bigcup_{n\in \N}A_n\in \mathcal{B}$.
%   \item If $A\in \mathcal{A}$ then $A^c\in \mathcal{A}$.
%\end{enumerate}

\medskip

Given a ring $\rr$, a \emph{premeasure} on it is a function $\mu :\rr \to [0, +\infty ]$ such that
\begin{enumerate}
    \item[($M_0$)] $\mu (\emptyset )=0\, $;
    \item[($M_1$)] for any sequence of pairwise disjoint sets $(A_n)_{n\in \N}$ in $\rr$ such that
    $\bigsqcup_{n\in \N} A_n \in \rr$, $$\mu \left( \bigsqcup_{n\in \N} A_n \right)=\sum_{n\in \N} \mu
    (A_n)\, .$$
\end{enumerate}

Property ($M_1$) is equivalent to
\begin{enumerate}
    \item[($M_1'$)] $\mu (A\sqcup B)=\mu (A) +\mu(B)$;
    \item[($M_1''$)]  If $(A_n)_{n\in\naturals}$ is a
non-increasing  sequence of sets in $\rr$ such that $\bigcap_{n\in \N} A_n=\emptyset$, then
$\lim_{n\to\infty}\mu(A_n)=0$.
\end{enumerate}

A premeasure is called $\sigma$-\emph{finite} if there exists a sequence $(A_n)$ in $\rr$ such that
$\mu (A_n)<+\infty$ for every $n$, and $\bigcup_{n} A_n =Y$.

A premeasure defined on a $\sigma$-algebra is called a \emph{measure}.

An {\em additive function} on a ring $\mathcal R$ is a map $\mu :\rr \to [0, +\infty ]$ satisfying properties $(M_0)$ and $(M_1)$.

%\begin{rem}[defining an additive function on a ring]\label{rem:additivefunction}
%Assume ${\mathcal A}\subset \rr$ is a subset of a ring such that $\mathcal A$ is stable by finite
%intersection, and each $A\in\rr$ can be written as a finite disjoint union $A=\sqcup_{i=1}^{n}A_i$,
%where each $A_i$ belongs to $\mathcal A$. Assume also that $\mu :\mathcal A \to [0, +\infty ]$  is a
%map such that for every decomposition $A=\sqcup_{i=1}^{i=n}A_i$ (where $A\in{\mathcal A}$ and each
%$A_i$ belongs to $\mathcal A$), we have $\mu(A)=\sum_{i=1}^{i=n}\mu(A_i)$. Then $\mu$ extends to a
%unique additive function on $\rr$.

%Indeed if $B\in\rr$ can be decomposed in two ways $B=\sqcup_{i=1}^{i=n}A_i$ and $B=\sqcup_{j=1}^{j=n}A'_j$, then we also have $B=\sqcup_{i,j}A_{ij}$, where we have set $A_{ij}=A_i\cap A'_j$. Since $\mathcal A$ is stable under intersection $A_{ij}\in\mathcal A$. Thus by the additive property of $\mu$ we have $\mu(A_i)=\sum_j\mu(A_{ij})$. Similarly  $\mu(A'_j)=\sum_i\mu(A_{ij})$, and thus $\sum_i\mu(A_i)=\sum_{i,j}\mu(A_{ij})=\sum_j\mu(A'_j)$.
%\end{rem}

We need a precise version of Caratheodory's Theorem on the extension of any premeasure $\mu$ to a measure,
therefore we recall here the notion of outer measure. For every $Q\subset Y$ let $\uu (Q)$ designate
the set of all sequences $(A_n)$ in $\rr$ such that $Q\subset \bigcup_{n} A_n$. Define
$\mu^*(Q)=+\infty $ if $\uu (Q)=\emptyset$; if  $\uu (Q)\neq \emptyset$ then
$$
\mu^* (Q)=\inf \left\{ \sum_{n=1}^\infty \mu (A_n) \; ;\; (A_n)\in \uu (Q) \right\}\, .
$$

The function $\mu^*$ is an \emph{outer measure} on the set $Y$.

A subset $A$ of $Y$ is called $\mu^*$-\emph{measurable} if for every $Q\in \pp (Y)$,
$$
\mu^* (Q)= \mu^* (Q\cap A)+\mu^* (Q\cap A^c)\, .
$$
\begin{thm}[Carath\'eodory \cite{HBauer:measure}, $\S I.5$]\label{thCar}
\begin{enumerate}
    \item The collection $\mathcal{A}^*$ of  $\mu^*$-measurable sets is a $\sigma$-algebra containing $\rr$, and
the restriction of $\mu^*$ to $\mathcal{A}^*$ is a measure, while the restriction of $\mu^*$ to $\rr$
coincides with $\mu$.
    \item If $\mu$ is $\sigma$-finite, then it has a unique extension to a measure on the
    $\sigma$-algebra generated by $\rr$.
\end{enumerate}
\end{thm}

\subsection{Spaces with measured walls.} From \cite{HaPa}, we recall that a \emph{wall} of a set $X$ is
 a partition $X=h\sqcup h^c$ (where {$h$ is possibly empty or the whole $X$}).  A collection $\hh$ of subsets of $X$
is called a \emph{collection of half-spaces} if for every $h\in \hh$ the complementary subset $h^c$ is
also in $\hh$. We call \emph{collection of walls} on $X$ the collection $\ww_\hh$ of pairs $w=\{h,h^c
\}$ with $h\in \hh$. For a wall $w=\{h,h^c\}$ we call $h$ and $h^c$ the two {\em half-spaces bounding
$w$}.

We say that a wall $w=\{h,h^c \}$ {\em separates} two disjoint subsets $A,B$ in $X$ if $A\subset h$ and
$B\subset h^c$ or vice-versa and denote by $\ww (A| B)$ the set of walls separating $A$ and $B$. In
particular $\ww (A |\emptyset)$ is the set of walls $w=\{h,h^c \}$ such that $A\subset h$ or $A\subset
h^c$; hence $\ww (\emptyset | \emptyset )=\ww$.

When $A=\{x_1,\dots,x_n\},B=\{y_1,\dots,y_m\}$ we write
$${\mathcal W}(A\vert B)={\mathcal
W}(x_1,\dots,x_n \vert y_1,\dots,y_m)\, .$$

We use the notation $\ww (x|y)$ to designate $\ww (\{x\}|\{y\})$. We call
 any set of walls of the form ${\mathcal W}(x\vert y)$ a {\em wall-interval}.
 By convention $\ww(A|A)=\emptyset$ for every non-empty set $A$.
\begin{defn}[space with measured walls \cite{CherixMartinValette}]\label{defn:spacemw}
A \emph{space with measured walls} is a 4-uple $(X, {\mathcal W}, {\mathcal B},\mu)$, where
$\mathcal{W}$ is a collection of walls, $\mathcal{B}$ is a $\sigma$-algebra of subsets in $\mathcal{W}$
and $\mu$ is a measure on $\mathcal{B}$, such that for every two points $x,y\in X$ the set of
separating walls $\ww (x| y)$ is in $\bb$ and it has finite measure. We denote by $\pdist_\mu$ the
pseudo-metric on $X$ defined by $\pdist_\mu (x,y)=\mu \left( \ww (x| y) \right)$, and we call it the
\emph{wall pseudo-metric}.
\end{defn}

\begin{lem}\label{lem:finite}
The collection $\rr$ of disjoint unions $\bigsqcup_{i=1}^n \ww (F_i|G_i)$, where $n\in \N^*$, and
$F_i,G_i$ are finite non-empty sets for every $i=1,2,...,n$, is a ring.
\end{lem}

\proof Property (1) is obviously satisfied.

We first note that $\ww (F|G) \cap \ww (F'|G') = \ww (F\cup F' | G\cup G') \sqcup \ww (F\cup G' | G\cup
F')$.

%Let now $F,G$ be two finite non-empty sets, and let $\mathcal{S}$ be the set of pairs $(S_1,S_2)$ such that $S_1 \subseteq F \, ,\, S_2 \subseteq G$ and $S_1\cup S_2 \subsetneq F\cup G$. It is easily seen that $$\ww(F|G)^c =\bigsqcup_{(S_1,S_2)\in \mathcal{S}} \ww (S_1 \cup S_2 | (F\cup G) \setminus (S_1 \cup S_2 ))\, .$$

Let now $F,G$ be two finite.

$$\ww(F|G)^c =\bigsqcup_{S\sqcup T = F\cup G,\{S,T\}\neq  \{F,G\}} \ww (S |  T )\, .$$

%and let $w=\{h,h^c\}$ be a wall not separating $F$ from $G$. Then either $w$ separates two points in $F$, or $w\in \ww(F\vert \emptyset)$ and $w$ separates   two points in $G$, or $\ww(F\cup G\vert \emptyset)$. We thus  have:
%$$\ww(F|G)^c =\bigl(\bigsqcup_{S_1\sqcup S_2 = F,S_i\neq\emptyset} \ww (S_1 |  S_2 )\bigr)\sqcup \bigl(\bigsqcup_{T_1\sqcup T_2 = G,T_i\neq\emptyset} \ww (F\cup T_1 |  T_2 )\bigr)\sqcup \ww(F\cup G\vert \emptyset)\, .$$

>From the two statements above it follows that $\rr$ satisfies property (3), i.e. it is closed with respect
to the operation $\setminus \, $. But $\rr$ is closed with respect to intersection,
thus $\rr$ is also closed with respect to union.\endproof

Theorem \ref{thCar} and Lemma \ref{lem:finite} imply the following.

\begin{prop}[minimal data required for a structure of measured walls]\label{prop:minimum}
Let $X$ be a space and let $\ww$ be a collection of walls on it. A structure of measured walls can be
defined on $(X,\ww )$ if and only if on the ring $\rr$ composed of disjoint unions $\bigsqcup_{i=1}^n
\ww (F_i|G_i)$, where $n\in \N^*$, and $F_i,G_i, i=1,2,...,n,$ are finite non-empty sets, can be
defined a premeasure $\mu $ such that for every $x,y\in X$, $\mu\left(\ww (x| y) \right)$ is finite.
\end{prop}

Let $(X,\ww , \bb ,\mu)$ and $(X',\ww' , \bb' ,\mu')$ be two spaces with measured walls, and let
$\phi:X\to X'$ be a map.
\begin{defn}\label{defn:hommesw}
The map $\phi$ is a \emph{homomorphism between spaces with measured walls} provided that:
\begin{itemize}
\item for any $w'=\{h',h'^c\}\in\ww'$ we have $\{\phi^{-1}(h'),\phi^{-1}(h'^c)\}\in\ww$ -
this latter wall we denote by $\phi^*(w')$;
\item the map $\phi^* : \ww' \to \ww$ is surjective and for every $B\in \bb$,
$(\phi^*)\iv (B) \in \bb'$ and $\mu' \left( (\phi^*)\iv (B) \right) = \mu (B)$.
\end{itemize}
Note that $\phi$ induces an isometry of the spaces equipped with the wall pseudo-distances.
\end{defn}
Consider the set $\hh$ of half-spaces determined by $ \ww$, and the natural projection map $\fp :\hh
\to \ww$, $h\mapsto \{h,h^c \}$. The pre-images of the sets in $\mathcal{B}$ define a $\sigma$-algebra
on $\hh$, {which we denote by $\bb^\hh$}; hence on $\hh$ can be defined a pull-back measure that we
also denote by $\mu$. This allows us to work either in $\hh$ or in $\mathcal{W}$. Notice that the $\sigma$-algebra $\bb^\hh$ does not separate points in $\hh$, as sets  in $\bb^\hh$ are unions of fibers of $\fp$.

\begin{defn}[\cite{ChatterjiNiblo}, \cite{Nica:cubulating}]
A section $\mathfrak{s}$ for $\fp$ is called \emph{admissible} if its image contains together with a
half-space $h$ all the half-spaces $h'$ containing $h$.
\end{defn}
Throughout the paper we identify an admissible section  $\fs$ with its image $\sigma=\fs (\ww )$; with
this identification, an admissible section becomes a collection of half-spaces, $\sigma$, such that:
\begin{itemize}
    \item  for every wall $w = \{h,h^c\}$ either $h$ or $h^c$ is in $\sigma $, but never both;
    \item  if $h\subset h'$ and $h\in \sigma$ then $h'\in \sigma$.
\end{itemize}

For any $x\in X$ we denote by $\fs_x$ the section of $\fp$ associating to each wall the half-space
bounding it and containing $x$. Obviously it is an admissible section. We denote by $\sigma_x$ its
image, that is the set of half-spaces $h\in\hh$ such that $x\in h$. Observe that $\sigma_x$ is not necessarily in $\bb^\hh$.

Note that $\fp(\sigma_x\vartriangle\sigma_y)=\ww (x| y)$.

\begin{exmp}[real hyperbolic space]\label{exmp:realhyper} For all the discussion below, see  \cite{Robertson:Crofton}.

Define the half-spaces of the  real hyperbolic space $\field H^n$ to be closed or open geometric
half-spaces, with boundary an isometric copy of $\field H^{n-1}$, so that a wall consists of one closed half-space and its (open) complement, as in Section 3 of \cite{CherixMartinValette}. Recall that the full group of direct
isometries of ${\field H^n}$ is $SO_0(n,1)$. The associated set of walls $\ww_{\field H^n}$ is
naturally identified with the homogeneous space $SO_0(n,1)/SO_0(n-1,1)$; as $SO_0(n-1,1)$ is
unimodular, there is a $SO_0(n,1)$--invariant borelian measure $\mu_{\field H^n}$ on the set of walls
\cite[Chapter 3, Corollary 4]{Nachbin:Haar}. The set of walls separating two points has compact closure
and finite measure. Thus $(\field H^n,\ww_{\field H^n},\bb,\mu_{\field H^n})$ is a space with measured
walls. By Crofton's formula \cite[Proposition 2.1]{Robertson:Crofton} up to multiplying the measure
$\mu_{\field H^n}$ by some positive constant the wall pseudo-metric on $\field H^n$ is just the usual
hyperbolic metric.
\end{exmp}

\begin{defn}\label{contactmw}
The action by automorphisms of a topological group $G$ on a space with measured walls  $(X,\ww , \bb
,\mu)$ is called \emph{continuous} if for every $x\in X$ the map $G \to X \, ,\, g\mapsto gx$ is
continuous, where $X$ is endowed with the topology defined by the pseudo-distance $\pdist_\mu$.
\end{defn}

The following result allows to produce many examples of spaces with measured walls.

\begin{lem}[pull back of a space with measured walls]\label{lem:pullback}
Let $(X,\ww , \bb ,\mu)$ be a space with measured walls, let $S$ be a set and $f:S\to X$ a map. There
exists a pull back structure of space with measured walls $(S, \ww_S ,\bb_S ,\mu_S)$ turning $f$ into a
homomorphism. Moreover:
\begin{itemize}
    \item[(i)] if $S$ is endowed with a pseudo-metric $\pdist$ and $f$ is an isometry between $(S,\pdist)$ and
$(X,\pdist_\mu)$, then the wall pseudo-metric $\pdist_{\mu_S}$ coincides with the initial pseudo-metric
$\pdist$;

\sm

    \item[(ii)] if a group $G$ acts on $S$ by bijective transformations and on $X$ by automorphisms of
    space with measured walls, and if $f$ is $G$-equivariant, then $G$ acts on $(S, \ww_S ,\bb_S
    ,\mu_S)$ by automorphisms of space with measured walls. Moreover, if the action on $X$ is continuous, the action on $S$ is.
\end{itemize}
\end{lem}

\proof Define the set of walls $\ww_S$ on $S$ as the set of walls $\{ f\iv (h), f\iv (h^c)\}$, where
$\{h , h^c\}$ is a wall in $X$. This defines a surjective map $f^* : \ww \to \ww_S$. We then consider
the push-forward structure of measured space on $\ww_S$. This defines a structure of measured space
with walls on $S$ such that $f$ is a homomorphism of spaces with measured walls.

(i) It is easily seen that for every $x,y\in S$, $(f^*)\iv (\ww_S (x|y))= \ww (f(x),f(y))$, hence
$\pdist_{\mu_S} (x,y)=\pdist_\mu (f(x),f(y))= \pdist (x,y)$.

(ii) If $f$ is $G$-equivariant then the whole structure of space with measured walls $(S, \ww_S ,\bb_S
    ,\mu_S)$ is $G$-equivariant.
\endproof

\medskip

One of the main interests in actions of groups on spaces with measured walls is given by the following
result.

\begin{lem}[\cite{CherixMartinValette}, \cite{CornTessValette:lp}]\label{lem:actlp}
Let $G$ be a group acting (continuously) by automorphisms on a space with measured walls $(X,\ww , \bb
, \mu )$. Let $p>0$ and let $\pi_p$ be the representation of $G$ on $L^p (\hh ,\mu_\hh )$.

Then for every $x\in X$, the map $b: G \to L^p (\hh ,\mu_\hh )$ defined by $b(g)=
\chi_{\sigma_{gx}}-\chi_{\sigma_{x}}$ is a (continuous) $1$-cocycle in $Z^1 (G, \pi_p)$. In other
words, a (continuous) action of $G$ on $L^p (\hh , \mu_\hh )$ by affine isometries can be defined by:
$$
g\cdot f = \pi_p(g)f + b(g)\, .
$$
\end{lem}

\begin{rem}
Recall that for a space $L^p (X, \mu )$ with $p\in (0,1)$, $\| f \|_p = \left( \int |f|^p d\mu
\right)^{\frac{1}{p}}$ no longer satisfies the usual triangular inequality, it only satisfies a similar
inequality with a multiplicative factor added to the second term. On the other hand, $\| f \|_p^p$ is
no longer a norm, but it does satisfy the triangular inequality, hence it defines a metric
\cite{KaltonPeck}.

In this paper we consider $L^p$-spaces endowed with this metric, for $p\in (0,1)$.
\end{rem}

\subsection{Embedding a space with measured walls in a median space.}\label{embmesmed}

Let $(X, {\mathcal W}, {\mathcal B},\mu)$ be a space with measured walls, and let $x_0$ be a base point
in $X$.

Recall from Example~\ref{exmp:med}, (\ref{exmp:symdiff}), that ${\bb^\hh}_{\sigma_{x_0}}$ denotes the
collection of subsets $A\subset \hh$ s.t. $A\vartriangle \sigma_{x_0}\in\bb$ and $\mu(A\vartriangle
\sigma_{x_0})<+\infty \, $, and that endowed with the pseudo-metric $\pdist_\mu(A,B)=\mu(A\vartriangle
B)$ this collection becomes a median pseudo-metric space. The map
\begin{equation}\label{chixo}
\chi^{x_0}:{\bb^\hh}_{\sigma_{x_0}}\to {\mathcal S}^1(\hh,\mu),\, \,  \chi^{x_0}
(A)=\chi_{A\vartriangle \sigma_{x_0}}
\end{equation}
 is an isometric embedding of ${\bb^\hh}_{\sigma_{x_0}}$ into the median subspace
${\mathcal S}^1(\hh,\mu)\subset {\mathcal L}^1(\hh,\mu)$, where ${\mathcal S}^1(\hh ,\mu)=\{\chi_B \;
;$ $B$ measurable and $\mu(B)<+\infty\}$.

The formula $A\vartriangle \sigma_{x_1}=(A\vartriangle \sigma_{x_0})\vartriangle
(\sigma_{x_0}\vartriangle \sigma_{x_1})$ and the fact that $\sigma_{x_0}\vartriangle \sigma_{x_1}$ is
measurable with finite measure shows that the median pseudo-metric spaces ${\bb^\hh}_{\sigma_{x_0}}$
and ${\bb^\hh}_{\sigma_{x_1}}$ are identical: we simply denote this space by ${\bb^\hh}_{X}$. In
particular $\sigma_x\in {\bb^\hh}_{X}$ for each $x\in X$.

For $x,y\in X$ we have $\pdist_\mu(x,y)=\mu(\sigma_x\vartriangle\sigma_y)$, thus $x\mapsto \sigma_x$ is
an isometric embedding of $X$ into $({\bb^\hh}_{X},\pdist_\mu)$. Composing with the isometry
$\chi^{x_0}: {\bb^\hh}_{X} \to {\mathcal S}^1(\hh,\mu)\, $, we get the following  well-known result
stating that a wall pseudo-distance is of type $1$, in the terminology of \cite[Troisi\`eme partie, $\S
2$]{BretDacunhaCastelKriv}:
\begin{lem}\label{lem:embedL1}
Let $(X,\ww , \bb ,\mu)$ be a space with measured walls, and $x_0\in X$ a base point in it. Then the
map $x\mapsto \chi_{\ww(x| x_0)}$ defines an isometry from $X$ to ${L}^1(\ww,\mu)$. Thus if the wall
pseudo-distance is a distance then $(X,\dist_\mu)$ is isometric to a subset of ${L}^1(\ww,\mu)$, and so
it is submedian.
\end{lem}

We could probably define the median space associated to a space with measured walls $(X,\ww , \bb
,\mu)$ to be the median hull of the isometric image of $X$ inside ${L}^1(\ww,\mu)$ (and then perhaps
take the closure in order to get a complete median space). We give here an alternative construction
which is more intrinsic.

\begin{notation}
We denote by $\overline{\mathcal M}(X)$ the set of admissible sections, and by ${\mathcal M}(X)$ the
intersection $\overline{\mathcal M}(X)\cap {\bb^\hh}_{X}$.

Every section $\sigma_x$ belongs to ${\mathcal M}(X)$, thus $X$ isometrically embeds in ${\mathcal
M}(X)$. We denote by $\iota:X\to\mm(X)$ this isometric embedding.
\end{notation}
\begin{prop}\label{prop:medianization}
 Let  $(X,\ww , \bb ,\mu)$ be
a space with measured walls.
\begin{enumerate}
    \item[(i)] The space ${\mathcal M}(X)$ defined as above is a median subspace of ${\bb^\hh}_{X}$.

    \sm

    \item[(ii)]  Any homomorphism $\phi :X \to X'$ between $X$ and
another space with measured walls $(X',\ww' , \bb' ,\mu')$ induces an isometry $\mm(X)\to\mm(X')$.

\sm

    \item[(iii)] In particular the group of automorphisms of $(X,\ww , \bb,\mu)$ acts by isometries on ${\mathcal
M}(X)$.
\end{enumerate}
\end{prop}
\begin{proof} (i)\quad Given an arbitrary triple $(\sigma_1,\sigma_2,\sigma_3)\in\mm(X)^3$, let us denote by $m(\sigma_1,\sigma_2,\sigma_3)$ the
set of half-spaces $h$ such that there exist at least two distinct indices $i,j\in\{1,2,3\}$ with
$h\in\sigma_i,h\in\sigma_j$. In other words $m(\sigma_1,\sigma_2,\sigma_3)=(\sigma_1\cap\sigma_2)\cup
(\sigma_1\cap\sigma_3)\cup (\sigma_2\cap\sigma_3)$ (see also Example \ref{boa}).

Clearly $m=m(\sigma_1,\sigma_2,\sigma_3)$ belongs to $\overline{\mathcal M}(X)$. Fix a point $x_0$ in
$X$ and take $\chi_0= \chi^{x_0}$ the function defined in (\ref{chixo}). We want to show that $\chi_0(m
)=m(\chi_0(\sigma_1),\chi_0 (\sigma_2),\chi_0 (\sigma_3))$. This will prove that $m\in{\bb^\hh}_{X}$
and that $m$ is a median point of $\sigma_1,\sigma_2,\sigma_3$.

For our set-theoretical calculation it is convenient to treat characteristic functions as maps from
$\hh$ to $\integers /2\integers$. We may then use the addition (mod. 2) and pointwise multiplication on
these functions. We get
$$\chi_{A\cap B}=\chi_A\chi_B,\: \chi_{A\vartriangle B} = \chi_A+\chi_B,\: \chi_{A\cup B}=\chi_A+\chi_B+\chi_A\chi_B\, .
$$
It follows easily that for any three subsets $A,B,C$ we
have
$$
\chi_{(A\cap B)\cup(A\cap C)\cup(B\cap C)} = \chi_A\chi_B+\chi_A\chi_C+\chi_B\chi_C\, . $$ Thus
$\chi_{[(A\cap B)\cup(A\cap C)\cup(B\cap C)]\vartriangle D} =
\chi_A\chi_B+\chi_A\chi_C+\chi_B\chi_C+\chi_D$. On the other hand $\chi_{((A\vartriangle D)\cap
(B\vartriangle D))\cup((A\vartriangle D)\cap (C\vartriangle D))\cup((B\vartriangle D)\cap
(C\vartriangle D))} =
(\chi_A+\chi_D)(\chi_B+\chi_D)+(\chi_A+\chi_D)(\chi_C+\chi_D)+(\chi_B+\chi_D)(\chi_C+\chi_D) =
\chi_A\chi_B+\chi_A\chi_C+\chi_B\chi_C +2\chi_A\chi_D+2\chi_B\chi_D+2\chi_C\chi_D+3\chi_D
=\chi_A\chi_B+\chi_A\chi_C+\chi_
 B\chi_C +\chi_D  $. We have thus checked that $[(A\cap B)\cup(A\cap C)\cup(B\cap C)]\vartriangle D$ coincides
 with $[(A\vartriangle D)\cap (B\vartriangle D)]\cup [(A\vartriangle D)\cap (C\vartriangle D)]\cup
 [(B\vartriangle D)\cap (C\vartriangle D)]$. Applying this to $A=\sigma_1,B=\sigma_2,C=\sigma_3,D=\sigma_{x_0}$ yields the desired result.

\medskip

(ii)\quad Consider a homomorphism of spaces with measured walls $\phi:X\to X'$. It is easily seen that
the surjective map $\phi^* : \ww' \to \ww$ induces a surjective map $\phi^* :\hh' \to \hh$ such that
for every $B\in \bb^\hh$, $(\phi^*)\iv (B)\in \bb^{\hh'}$ and $\mu' \left( (\phi^*)\iv
(B)\right)=\mu(B)$.

Let $\sigma$ denote any admissible section. Set $\phi_*(\sigma)=\left( \phi^* \right)\iv (\sigma)=
\{h'\in \hh'\; ;\; \phi\iv (h')\in \sigma \}$. Since $\phi$ is a homomorphism, $\phi_*(\sigma)$ is an
admissible section of $(X',\ww' , \bb' ,\mu')$. Note that $\phi_*(\sigma_x)=\sigma_{\phi(x)}$ and that
$\phi_*(\sigma \vartriangle \sigma ')= \phi_*(\sigma ) \vartriangle \phi_*(\sigma')$. This implies that
$\phi_*$ defines a map from $\mm (X)$ to $\mm (X')$. Moreover
$\pdist_{\mm(X')}(\phi_*(\sigma),\phi_*(\sigma'))= \mu' (\phi_*(\sigma ) \vartriangle
\phi_*(\sigma'))=\mu'(\phi_*(\sigma \vartriangle \sigma '))= \mu'(\left( \phi^* \right)\iv (\sigma
\vartriangle \sigma '))= \mu (\sigma \vartriangle \sigma ')=\pdist_{\mm(X)}(\sigma,\sigma')$. Thus
$\phi_*$ is an isometry.

The statement (iii) is an immediate consequence of (ii).\end{proof}
The results in Proposition \ref{prop:medianization} justify the following terminology.
\begin{defn}
We call ${\mathcal M}(X)$  \emph{the median space associated with $(X,\ww , \bb ,\mu)$}.
\end{defn}

The first part of Theorem~\ref{res} is proved.

\begin{rem}\label{rem:mwallm(X)}
The median space $\mm(X)$ has measured walls. Indeed for each $h\in\hh$ define $h_\mm$ to be the set of
$\sigma\in\mm(X)$ such that $h\in\sigma$. The complement of $h_\mm$ in $\mm(X)$ is the set of
$\sigma\in\mm(X)$ such that $h\not\in\sigma$, or equivalently by the properties of admissible sections
$h^c\in\sigma$. In other words $({h_\mm})^c=(h^c)_\mm$. Thus $\{h_\mm\}_{h\in\hh}$ is a collection of
half-spaces - which we will denote by $\hh_\mm$. We denote by $\ww_\hh$ the associated set of walls on
$\mm(X)$. Using the bijection $\ww\to\ww_\hh$ induced by $h\mapsto h_\mm$ we define on $\ww_\hh$ a
$\sigma$-algebra $\bb_\hh$ and a measure $\mu_\mm$. Note that $\iota:X\to\mm(X)$ is a homomorphism.
Note also that the distance on $\mm (X)$ coincides with the distance induced by the measured walls
structure.

It is easy to check that the medianized space associated with $\mm(X)$ endowed with this structure of
space with measured walls is $\mm(X)$ itself.
\end{rem}

\begin{rem}\label{rembasic}

One cannot hope to define a median space $(\mm (X),\dist )$ associated to a space with measured walls
such that there exists an isometric map $\iota :(X, \pdist_\mu)\to (\mm (X),\dist )$ with the
universality property that any isometric map from $(X, \pdist_\mu)$ to a median space factors through
$\iota \, $. This was explained in Remark \ref{nofunctor}.
\end{rem}

%
%%%%%%%%%%%%%%%%%%%%%%%%%%%%%%%%%%%%%%%%%%%%%%%%%%%

%%%%%%%%%%%%%%%%%%%%%%%%%%%%%%%%%%%%%%%%%%%%%%%%
\section{A review of median algebras.}\label{sec:medianalgebras}
%%%%%%%%%%%%%%%%%%%%%%%%%%%%%%%%%%%%%%%%%%%%%%%%%%%%%%%%%%
The notion of median algebra appeared as a common generalization of trees and  lattices (in the ordered structure sense of the word). We recall here some basic definitions and properties related to median algebras. For proofs and further details we
refer the reader to the books \cite{VandeVel:book}, \cite{Verheul:book},                             the surveys \cite{BandeltHedlikova}, \cite{Isbell},
as well as the papers \cite{BirkhoffKiss}, \cite{Sholander1}, \cite{Sholander2} and \cite{Roller:median}.

\subsection{Definitions, examples.}
\begin{defn}(median algebra, first definition)
A \emph{median algebra} is a set $X$ endowed with a ternary operation $(a,b,c)\mapsto m(a,b,c)$ such
that:\begin{itemize}
    \item[(1)] $m(a,a,b)=a$;
    \item[(2)] $m(a,b,c)=m(b,a,c)=m(b,c,a)$;
    \item[(3)] $m(m(a,b,c),d,e)= m(a,m(b,d,e), m(c,d,e))$.
\end{itemize}
Property (3) can be replaced by $(3')$ $m(a,m(a,c,d),m(b,c,d))=m(a,c,d)$.

The element $m(a,b,c)$ is the \emph{median of the points }$a,b,c$. In a median algebra $(X,m)$, given any two points $a,b$ the set $I(a,b)=\{x\; ;\; x=m(a,b,x)\}$ is called \emph{the interval of endpoints} $a,b$. This defines a map $I :X\times X \to \pp (X)$. We say
that a point $x\in I(a,b)$ is \emph{between $a$ and $b$}.

A \emph{homomorphism} of median algebras is a map $f: (X,m_X)\to (Y,m_Y)$ such that
$m_Y(f(x),f(y),f(z))=f(m_X(x,y,z))$. Equivalently, $f$ is a homomorphism if and only if it preserves
the betweenness relation.
If moreover $f$ is injective (bijective) then $f$ is called \emph{embedding} or \emph{monomorphism}
(respectively \emph{isomorphism}) of median algebras.

\end{defn}

The following are straightforward properties that can be found in the
literature (see for instance \cite{Sholander1} and \cite[$\S 2$]{Roller:median}).

\begin{lem}\label{lp1} Let $(X,m)$ be a median algebra. For $x,y,z\in X$ we have that
\begin{enumerate}
    \item $I(x,x)=\{x\}$;
    \item $I(x,y)\cap I(x,z)=I(x,m(x,y,z))$;
    \item $I(x,y)\cap I(x,z) \cap I(y,z) = \{m(x,y,z)\}$;
    \item if $a\in I(x,y)$ then for any $t$, $I(x,t)\cap I(y,t)\subseteq I(a,t)$ (equivalently $m(x,y,t)\in
    I(a,t)\, $);
    \item if $x\in I(a,b)$ and $y\in I(x,b)$ then $x\in I(a,y)$.
\end{enumerate}
\end{lem}
A sequence of points $(a_1,a_2,...,a_n)$ is \emph{geodesic} in the median algebra $(X,m)$ if $a_i\in
I(a_{1},a_{i+1})$ for all $i=2,\dots,n-1$. This is equivalent, by Lemma \ref{lp1}, point (5), to the
condition that $a_{i+1}\in I(a_i, a_n)$ for all $i=1,2,...,n-2$.
\begin{lem}
If $(x,t,y)$ is a geodesic sequence, then:
\begin{enumerate}
\item $I(x,t)\cup I(t,y)\subseteq I(x,y)$;
\item $I(x,t)\cap I(t,y) = \{ t\}$.
\end{enumerate}
\end{lem}
%%%%%%%%%%
According to \cite{Sholander1}, \cite{Sholander2} there is an alternative definition of median
algebras, using intervals.

\begin{defn}(median algebra, second definition)\label{defn:equivmedian} A \emph{median algebra} is a set $X$ endowed with a map $I:X\times X \to
\pp (X)$ such that:
\begin{itemize}
    \item[(1)] $I(x,x)=\{ x\}$;
    \item[(2)] if $y\in I(x,z)$ then $I(x,y)\subset I(x,z)$;
    \item[(3)] for every $x,y,z$ in $X$ the intersection $I(x,y)\cap I(x,z)\cap I(y,z)$ has cardinality
$1$.
\end{itemize}
\end{defn}

\begin{exmp}\label{exmp:mmetricmalg}
Let $(X,\dist)$ be a median space. Then the metric intervals $I(x,y)$ satisfy the properties in
Definition \ref{defn:equivmedian}, and thus the metric median $(x,y,z)\mapsto m(x,y,z)$ defines a
structure of  median algebra on $X$.

\end{exmp}

\begin{exmp}\label{boa} Here is the set-theoretic generalization of Example~\ref{exmp:med}, (\ref{exmp:symdiff}).
For any set $X$, the power set $\pp (X)$ is a median algebra when endowed with the Boolean median
operation
\begin{equation}\label{boolemedian}
m(A,B,C)= (A\cap B) \cup (A\cap C) \cup (B\cap C) = (A\cup B) \cap (A\cup C) \cap (B\cup C)\, .
\end{equation}
The median algebra $(\pp (X), m)$ is called a \emph{Boolean median algebra}. One easily sees that in this case
\begin{equation}\label{booleaninterval}
I(A,B)= \{ C\; ;\; A\cap B \subset C \subset A\cup B\}\, .
\end{equation}
In what follows we use the notation $Bm(A,B,C)$ to designate the Boolean median defined in
(\ref{boolemedian}) and $BI(A,B)$ to designate the Boolean interval defined in (\ref{booleaninterval}).
\end{exmp}
It appears that Example \ref{boa} is in some sense the typical example of median algebra. More
precisely, according to Corollary \ref{sigmaemb}, any median algebra is a subalgebra of a  Boolean
median algebra, up to isomorphism.

\subsection{Convexity.}

\begin{defn}A \emph{convex subset} $A$ in a median algebra is a subset such that for any $a,b\in A $,
$I(a,b)\subset A$; equivalently it is a subset such that for every $x\in X$, and $a,b$ in $A$ the
element $m(a,x,b)$ is in $A$.

A subset $h$ in a median space $(X,m)$ is called \emph{a convex half-space} if itself and the
complementary set $h^c$ are convex. The pair $\{h, h^c\}$ is called a \emph{convex wall}. We denote by
$\hh_c(X)$ the set of convex half-spaces in $X$ and by $\ww_c(X)$ the set of convex walls in $X$. When there is no possibility of confusion we simply use the notations $\hh_c$ and $\ww_c$.
\end{defn}

The above algebraic notion of convexity coincides with the metric notion of convexity introduced in
Definition~\ref{conv}, in the case of the median algebra associated with a median space (see
Example~\ref{exmp:mmetricmalg}).

The following result shows that there are plenty of convex walls in a median algebra.
\begin{thm}\label{tsepp}
Let $X$ be a median algebra, and let $A, B$ be two convex non-empty disjoint subsets of $X$. Then there
exists a convex wall separating $A$ and $B$.\end{thm}
A proof of Theorem \ref{tsepp} when $A$ is a singleton can be found in \cite{Nieminen}; in its most
general form it follows from \cite[Theorem 2.5]{VandeVel}. Other proofs can be found in \cite[$\S
5.2$]{Basarab} and in \cite[$\S 2$]{Roller:median}.
\begin{cor}\label{cor:convexwalls}
Given any two distinct points $x,y$ in a median space $(X,\dist)$ there exists a convex wall
$w=\{h,h^c\}$ with $x\in h,y\in h^c$.
\end{cor}
\begin{defn}\label{defsigma1}
Given a median algebra $X$, one can define the map
$$\sigma :X \to \pp (\hh_c) \, ,\, \sigma(x)= \sigma_x=\{ h\in \hh_c\; ;\; x\in h\}.$$ \end{defn}
A
consequence of Theorem \ref{tsepp} is the following.
\begin{cor}\label{sigmaemb}
The map $\sigma$ is an embedding of median algebras.
\end{cor}

%%%%%%%%%%%%%%%%%%%%%%%%%%
\section{Median spaces have measured walls.}\label{section:medianmeswalls}\label{med=mws}
%%%%%%%%%%%%%%%%%%%%%%%%%%

%%%%%%%%%%%%%%%%%%%%%%%%%%%%%%%%%%%%%%%%%%%%%%%%%%%%%%%%%%%%%%%%%%%%%%%%%%%%%%%%%%%%%%%%%%

%%%%%%%%%%%%%%%%%%%%%%%%%%%%%%%%%%%%%%%%%%%%%%%%%%%%%%%%%%%%%%%%%%%%%%%%%%%%%%%%%%%%%%%%%%%%%
The aim of this section is to prove the following.
\begin{thm}\label{medwalls}
Let $(X,\dist)$ be a median space. Let $\ww$ be the set of convex walls, and let $\bb$ be the
$\sigma$-algebra generated by the following subset of $\pp (\ww )$:
$$\uu=\{ \ww(x|y)\; ;\; x,y\hbox{ points of }X\}\, .
$$
Then there exists a measure $\mu$ on ${\mathcal B}$ such that:
\begin{enumerate}
    \item  $\mu (\ww(x|y))=\dist (x,y)$;
consequently the 4-tuple $(X,\ww,{\mathcal B},\mu)$ is a space with measured walls;
    \item any isometry of $(X, \dist )$ is an automorphism of the space with measured walls $(X,\ww,{\mathcal
    B},\mu)$.
\end{enumerate}
\end{thm}

\begin{rem}\label{remunicity}
According to Caratheodory's theorem, a measure $\mu$ on the $\sigma$-algebra ${\mathcal B}$ is not
uniquely defined by the condition (1) in Theorem \ref{medwalls}. It is uniquely defined if there exists
say a sequence of points $(x_n)$ in $X$ such that $\ww = \bigcup_{n,m}\ww (x_n|x_m)$. This happens for
instance if there exists a countable subset in $X$ whose convex hull is the entire $X$. Uniqueness is also
guaranteed when for some topology on $\ww$ the measure $\mu $ is borelian and $\ww$ is locally compact
second countable.
\end{rem}
Combining Theorem \ref{medwalls} above and Lemma~\ref{lem:embedL1} we get the following:
\begin{cor}\label{cor:medinL1}
Let $(X,\dist)$ be a median space. Then $X$ isometrically embeds in $L^1(\ww , \mu )$, where $(\ww ,
\mu )$ are as in Theorem \ref{medwalls}.

More precisely, given any $x_0\in X$, the space $X$ is isometric to $\left\{\chi_{\ww(x|x_0)}\; ;\;
x\in X \right\}\subset L^1(\ww,\mu)$ endowed with the induced metric.
\end{cor}
 {The fact that median spaces embed isometrically into $L^1$--spaces was known previously, though not via a construction of an embedding as above, but using Assouad's result that a space is embeddable into an $L^1$--space if and only if any finite subset of the space is (\cite{Assouad:analytique}, \cite{Assouad:analytiqueComment}, \cite{AssouadDeza}). That finite median spaces can be embedded into $\ell^1$--spaces seems to be well known in graph theory; all proofs usually refer to finite median graphs only, but can be adapted to work for finite median spaces (see for instance \cite{Mulder:book}). There exist even algorithms which isometrically embed a given median graph into an $\ell^1$--space; the same method yields algorithms in sub-quadratic time recognizing median graphs \cite{Hagauer}. The statement that finite median spaces can be embedded into $\ell^1$ was explicitly stated and proved for the first time in \cite[Theorem V.2.3]{Verheul:book}.}
\begin{cor}\label{cor:premedian}
A metric space $(X,\dist )$ is submedian in the sense of Definition \ref{defn:submed} if and only if it
admits a structure of space with measured walls $(X, \ww ,\bb, \mu )$ such that $\dist = \dist_\mu$.
Moreover all walls in $\ww$ may be assumed to be convex.
\end{cor}

\proof The direct part follows from Theorem \ref{medwalls} and Lemma \ref{lem:pullback}.

The converse part follows from Lemma \ref{lem:embedL1}.\endproof
\begin{rem}\label{finite}
Corollary \ref{cor:premedian} for finite metric spaces was already known. We recall this
version here as it will prove useful further on.

More precisely, according to \cite{Assouad:analytique} and \cite{AssouadDeza} a finite metric space
$(X, \dist )$ is isometrically $\ell^1$-embeddable if and only if
$$\dist = \sum_{S\subseteq X}\lambda_S \delta_S \, ,\, $$
where $\lambda_S$ are non-negative real numbers, and $\delta_S(x,y)=1$ if $x\neq y$ and $S\cap \{x,y\}$
has cardinality one,  $\delta_S(x,y)=0$ otherwise.
\end{rem}
Theorem~\ref{medwalls} together with Proposition~\ref{prop:medianization} show that  the natural dual
category of median pseudo-metric spaces is the category of spaces with measured walls. Precise results
on duality of categories for particular categories of median algebras and spaces with walls can be
found in \cite{Roller:median} and in \cite{Basarab}.

\begin{rem}
According to the construction in $\S 3.2$, a space with measured walls $X$ has a natural embedding into
a median space $\mathcal{M} (X)$; moreover  $\mathcal{M} (X)$ has an induced structure of space with
measured walls, and its metric coincides with the metric induced by the measured walls structure
(Remark \ref{rem:mwallm(X)}).

We note here that the above structure of space with measured walls on  $\mathcal{M} (X)$ does not in
general agree with the structure described in this section. In general the first structure does not
have convex walls, as the walls on $X$ may not be convex. In a forthcoming paper we will show that both structures on $\mathcal{M} (X)$ are equivalent, in the sense that they induce the same structure of measured spaces with walls on finite subsets.
\end{rem}

\medskip
The strategy of the proof of Theorem~\ref{medwalls} is to use Proposition \ref{prop:minimum}. We first
show that for any pair of finite non-empty sets $F,G$ in $X$, $\ww (F|G)$ is equal to $\ww (a|b)$ for
some pair of points $a,b$. In order to do this we need the following intermediate results.
\begin{lem}\label{lem:wallgeod} Let $(x,y,z)$ be a geodesic sequence. Then we have the following
decomposition as a disjoint union:
$${\mathcal W}(x\vert z)={\mathcal W}(x\vert y)\sqcup {\mathcal W}(y\vert z).$$
\end{lem}
\begin{proof}First notice that by convexity of half-spaces, the intersection ${\mathcal W}(x\vert y)\cap {\mathcal W}(y\vert z)$ is empty.
Then the inclusion ${\mathcal W}(x\vert z)\subseteq {\mathcal W}(x\vert y)\cup{\mathcal W}(y\vert z)$ is clear because if a half-space $h$ contains $x$ but does not contain $z$, then either $h$ contains $y$ (in which case the wall $\{h,h^c\}$ separates $y$ from $z$) or $h^c$ contains $y$ (in
which case the wall $\{h,h^c\}$ separates $x$ from $y$). The inclusion ${\mathcal
W}(x\vert y)\cup {\mathcal W}(y\vert z)\subseteq {\mathcal W}(x\vert z)$ holds because if $h$ contains $x$ and $y\not\in h$, again by convexity we
cannot have $z\in h$ and hence $\{h,h^c\}$ separates $x$ from $z$. \end{proof}

As an immediate consequence we get the following:

\begin{cor}\label{cor:wallgeod}
For any geodesic sequence $(x_1,x_2,...,x_n)$ we have the following decomposition:
$${\mathcal W}(x_1\vert x_n)={\mathcal W}(x_1\vert x_2)\sqcup \cdots  \sqcup {\mathcal W}(x_{n-1}\vert x_n).$$
\end{cor}

\begin{cor}\label{cor:wquad}
If $(x,y)$ and $(x',y')$ are parallel pairs then $${\mathcal W}(x\vert y)={\mathcal W}(x'\vert
y')={\mathcal W}(x,x'\vert y,y')\, .$$
and $${\mathcal W}(x\vert y')={\mathcal W}(x'\vert
y)={\mathcal W}(x\vert y)\sqcup {\mathcal W}(x\vert x')\, .$$

\end{cor}
\begin{lem}\label{lem:w12}
Given three points $x,y,z$ with median point $m$, we have
${\mathcal W}(x\vert y,z)={\mathcal W}(x\vert m)$.
\end{lem}
\begin{proof}According to Lemma~\ref{lem:wallgeod} we have that ${\mathcal W}(x\vert
y)={\mathcal W}(x\vert m)\sqcup {\mathcal W}(m\vert y)$ and that
${\mathcal W}(x\vert z)={\mathcal W}(x\vert m)\sqcup {\mathcal
W}(m\vert z)$. It follows that
$${\mathcal W}(x\vert y,z)={\mathcal
W}(x\vert y)\cap {\mathcal W}(x\vert z)={\mathcal W}(x\vert
m)\sqcup ({\mathcal W}(m\vert y)\cap {\mathcal W}(m\vert z)).$$ But by convexity of the walls ${\mathcal
W}(m\vert y)\cap {\mathcal W}(m\vert z)=\emptyset$, and we are
done.\end{proof}

We will use intensively the following two operations:
\begin{defn}[projection and straightening]\label{defn:procedures}
Let $(x,y),(a,b)$ be two pairs of points of a median space $X$.
%We call the first pair $(x,y)$ the {\em input pair}, and we call the last pair $(a,b)$ the {\em target pair}.

The {\em  projection of $(x,y)$ with target $(a,b)$} is the pair $(x',y')$ defined by $x'=m(x,a,b),y'=m(y,a,b)$.

If furthermore $x,y\in I(a,b)$ we also consider the {\em straightening of the path $(a,x,y,b)$}, which  by definition is the path $(a,p,q,b)$, where  the pair $(p,q)$ is defined by $p=m(a,x,y),q=m(b,x,y)$.

\end{defn}
Observe that given two pairs of points $(x,y),(a,b)$, the central rectangle $[x',a',y',b']$ associated with $[x,a,y,b]$ (as defined in
Definition \ref{defn:centrect}) is obtained by first  projecting  $(x,y)$ with target $(a,b)$ - this yields the pair $(x',y')$ - and then  straightening  $(a,x',y,',b)$ - which yields the pair $(a',b')$. We now give some properties of both procedures.

\begin{lem}\label{lem:wallproject}
Let $(x,y),(a,b)$ be two pairs of points.
\begin{enumerate}
\item Let $(x',y')$ be the projection of $(x,y)$ with target $(a,b)$. Then
$$\ww (x'|y')={\mathcal W}(x\vert y)\cap {\mathcal W}(a\vert b) \, .$$
\item Assume $x,y\in I(a,b)$, and let $(p,q)$ be the projection of $(a,b)$ with target $(x,y)$. Then $[p,x,q,y]$ is a rectangle, $\ww(p\vert q)=\ww(x\vert y)$, and $(a,p,q,b)$ is a geodesic sequence (thus $(a,x,y,b)$ has really been straightened to a geodesic).
\item\label{propertiesrectangle} Let $[x',a',y',b']$ be the central rectangle associated with $[x,a,y,b]$. Then
$$\ww (x'|y')={\mathcal W}(x\vert y)\cap {\mathcal W}(a\vert b), \ww (x'|y')=\ww (a'|b') \, .$$
\end{enumerate}
\end{lem}
%keep these notations because of the application
\begin{proof}
Since the central rectangle is in fact obtained by composing the projecting and straightening
operations, it is enough to prove statement~\ref{propertiesrectangle}.

The equality $\ww (x'|y')=\ww (a'|b')$ follows by  Corollary \ref{cor:wquad}.

By Lemma~\ref{lem:w12} we have ${\mathcal W}(x\vert x')={\mathcal W}(x\vert a,b)$. In particular
${\mathcal W}(x\vert x')\cap {\mathcal W}(a\vert b)=\emptyset$. And similarly ${\mathcal W}(y\vert
y')\cap {\mathcal W}(a\vert b)=\emptyset$.

Consider now a half-space $h$ such that $x\in h,y\not\in h$ and $\{h,h^c\}\in {\mathcal W}(a\vert b) $.
Since ${\mathcal W}(x\vert x')\cap {\mathcal W}(a\vert b)=\emptyset$, we deduce that $x'\in h$.
Similarly we have $y'\in h^c$. We have thus proved that ${\mathcal W}(x\vert y)\cap {\mathcal W}(a\vert b)\subset {\mathcal W}(x'\vert y')$.

On the other hand, since $\ww (x'|y') = \ww (a'|b')$ and $(a,a',b',b)$ is a geodesic, it follows that
$\ww (x'|y') \subset \ww (a|b)$.

According to Lemma \ref{RectCan2}, $(x',y')$ is parallel to a pair $(x'',y'')$ such that
$(x,x'',y'',y)$ is geodesic. This and Corollary \ref{cor:wquad} imply that $\ww (x'|y') \subset \ww
(x|y)$.\end{proof}

\begin{prop}\label{prop:complementstable}
Let $F$ and $G$ be two finite non-empty subsets in $X$. There exist two points $p,q\in X$ such that
$${\mathcal W}(F|G) = {\mathcal W}(p|q)\, .$$
\end{prop}
\begin{proof}
We use an inductive argument over $n=\card F + \card G$. For $n=2$ the result is obvious, while for
$n=3$ it is Lemma \ref{lem:w12}.

Assume that the statement holds for $n$ and let $F,G$ be such that $\card F + \card G=n+1\geq 3$.
Without loss of generality we may assume that $\card F \geq 2$. Then $F= F_1 \sqcup \{x\}$, and $\ww
(F|G)= \ww (F_1 |G) \cap \ww (x|G)$. The inductive hypothesis implies that $\ww (F_1 |G)= \ww (a|b)$
and $\ww (x|G)= \ww (c|d)$, for some points $a,b,c,d$. Hence $\ww (F|G)= \ww (a|b) \cap \ww (c|d)$. We
end up by applying Lemma \ref{lem:wallproject}.\end{proof}
At this stage we have proven that the ring $\rr$ defined in  Proposition \ref{prop:minimum} coincides
with the set of disjoint unions $\bigsqcup_{i=1}^n \ww (x_i|y_i)$.  It remains to show that there is a premeasure $\mu:\rr\to\R^+$ on the ring $\rr$ such that $\mu(\ww(x\vert y))=\dist(x,y)$. We first
define $\mu$ as an additive function.

\begin{lem}\label{lem:induct}
If $\ww (x|y)=\ww (a|b)$ then $\dist (x,y)= \dist (a,b)$.
\end{lem}

\proof First let $(x',y')$ be the projection of $(x,y)$ with target $(a,b)$. Then by Lemma \ref{lem:wallproject}(1) we have $\ww (x'|y') = {\mathcal W}(x\vert y)\cap {\mathcal W}(a\vert b) ={\mathcal W}(a\vert b)$. By Corollary~\ref{cor:medlip} the median map is 1-Lipschitz, thus $d(x',y')\le d(x,y)$.

We now straighten $(a,x',y',b)$ to $(a,p,q,b)$ (thus $(p,q)$ is the projection of $(a,b)$ with target $(x',y')$). Then by Lemma \ref{lem:wallproject}(2) we have $\ww(p \vert q)=\ww(x' \vert y')=\ww(a\vert b)$, and $(a,p,q,b)$ is a geodesic.  By Corollary~\ref{cor:wallgeod} we deduce $\ww (a | p)=\ww (q|b)=\emptyset$, and thus $a=p,q=b$. It follows that $d(a,b)=d(p,q)$, and thus by Corollary~\ref{cor:wquad} we have $d(a,b)=d(x',y')\le d(x,y)$. We conclude by symmetry.
\endproof
\begin{prop}\label{muwelldefined} Assume that for two points $a,b$ the set of walls $\ww (a|b)$
decomposes as ${\mathcal W}(a\vert b)=\bigsqcup_{j=1}^{n} {\mathcal W}(x_j|y_j)$. Then there exists a
geodesic sequence $(a_1=a,a_2,\dots,a_{2^n}=b)$ and a partition $\{1,2,\dots,2^n-1\}=I_1\sqcup
I_2\sqcup\dots\sqcup I_n$ such that:
\begin{enumerate}
\item for each $j\in\{1,\dots,n\}$ the set $I_j$ has $2^{j-1}$ elements and we have a decomposition of
${\mathcal W}(x_j\vert y_j)= \bigsqcup _{i\in I_j}{\mathcal W}(a_{i} \vert  a_{i+1})$
\item for each $j\in\{1,\dots,n\}$ we have $\dist(x_j,y_j)=\sum_{i\in I_j}\dist(a_{i},a_{i+1})$
\end{enumerate}
In particular, $\dist(a,b)=\sum_j \dist(x_j,y_j)$.
\end{prop}
We easily deduce the following:
\begin{cor}
There is a unique additive function $\mu:\rr \to\R^+ $ such that $\mu (\ww(x|y))=\dist(x,y)$.
\end{cor}

\begin{figure}[!ht]
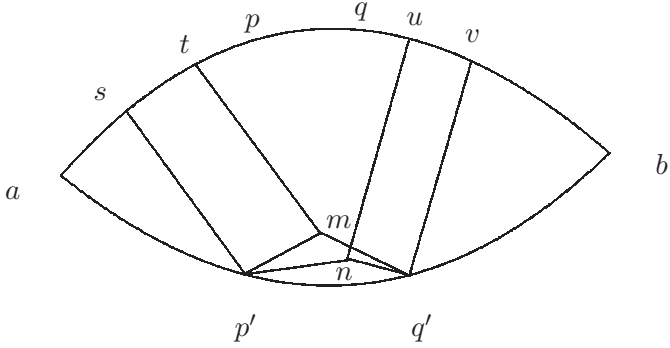

\centering \caption{The construction in Lemma \ref{aux}.} \label{fig3}
\end{figure}

To prove the Proposition we need the following auxiliary result:
\begin{lem}\label{aux} In a median space $(X,\dist)$, consider two geodesic sequences with common
endpoints $(a,p,q,b)$ and $(a,p',q',b)$, such that ${\mathcal W}(p\vert q)\cap {\mathcal W}(p'\vert
q')=\emptyset$. Let $(s,t)$ be the projection of $(p',q')$ with target $(a,p)$. Similarly let  $(u,v)$ be the projection of $(p',q')$ with target $(q,b)$. Then $\dist(p',q')=\dist(s,t)+\dist(u,v)$.\end{lem}
\begin{proof} Consider two more points: $m=m(t,p',q'),
n=m(u,p',q')$ (see Figure \ref{fig3}). Let us check that $[s,t,m,p']$ is a rectangle. By construction
$(t,m,p')$ is a geodesic sequence. Since $s,t$ are projection of $p',q'$ onto the interval $I(a,p)$ we deduce that $(q',m,t,s),(p',s,t)$ are geodesic sequences. And since $(x,p',q',y)$ is a geodesic sequence we see that $(x,s,p',m,q',y)$ is geodesic. %Similarly one shows that $[u,v,q',n]$ is a rectangle.

We thus have $\dist(p',m)=\dist(s,t)$, and also $\ww(p'\vert m)=\ww(s\vert t)$ (by Corollary~\ref{cor:wquad}). Hence $\ww(p'\vert m)=\ww(a\vert p)\cap\ww(p'\vert q')$ (by Lemma~\ref{lem:wallproject}(1)). Similarly we get $\dist(n,q')=\dist(u,v)$, and $\ww(n\vert q')=\ww(q\vert b)\cap \ww(p'\vert q')$.

 We claim that $\ww(m\vert q')=\ww(q\vert b)\cap \ww(p'\vert q')$. Indeed applying several times Lemma~\ref{lem:wallgeod} we get
$$\ww(p'\vert m)\sqcup \ww(m\vert q')=\ww(p'\vert q')\subset \ww(a\vert b)=\ww(a\vert p)\sqcup \ww(p\vert q)\sqcup \ww(q\vert b)$$ and the claim follows, since by assumption $\ww(p\vert q)\cap \ww(p'\vert q')=\emptyset$ and we already have $\ww(p'\vert m)=\ww(a\vert p)\cap\ww(p'\vert q')$.

We deduce that $\ww(m\vert q')=\ww(n\vert q')$. This implies $\dist(m,q')=\dist(n,q')=\dist(u,v)$ by Lemma~\ref{lem:induct}. Since $(p',m,q')$ is a geodesic we get $\dist(p',q')=\dist(p',m)+\dist(m,q')=\dist(s,t)+\dist(u,v)$.\end{proof}

%TeXCAD Options
%\grade{\on}
%\emlines{\off}
%\epic{\off}
%\beziermacro{\on}
%\reduce{\on}
%\snapping{\off}
%\quality{8.00}
%\graddiff{0.01}
%\snapasp{1}
%\zoom{4.0000}
\unitlength 0.75mm % = 2.85pt
\linethickness{0.4pt}
\ifx\plotpoint\undefined\newsavebox{\plotpoint}\fi % GNUPLOT compatibility
\begin{picture}(130.25,90)(0,0)
\qbezier(23.75,45)(67.63,95)(121,49) \qbezier(121,49)(74.63,4)(23.75,45)
\put(15.25,42){\makebox(0,0)[cc]{$a$}} \put(130.25,47){\makebox(0,0)[cc]{$b$}}
\put(57.75,72){\makebox(0,0)[cc]{$p$}} \put(77,74){\makebox(0,0)[cc]{$q$}}
\put(56.5,18){\makebox(0,0)[cc]{$p'$}} \put(87.75,18){\makebox(0,0)[cc]{$q'$}}
%\emline(35.5,96.25)(56.25,67.75)
\multiput(35.5,56.25)(.033739837,-.046341463){615}{\line(0,-1){.046341463}}
%\end
%\emline(56.25,67.75)(56.25,67.75)
\put(56.25,27.75){\line(0,1){0}}
%\end
%\emline(85.75,67.5)(85.75,68)
\put(85.75,27.5){\line(0,1){.5}}
%\end
\put(30.75,59.5){\makebox(0,0)[cc]{$s$}} \put(45.75,68.25){\makebox(0,0)[cc]{$t$}}
\put(86.5,72.75){\makebox(0,0)[cc]{$u$}}
%\emline(47.75,104.5)(69.75,74.75)
\multiput(47.75,64.5)(.033690658,-.045558959){653}{\line(0,-1){.045558959}}
%\end
%\emline(69.75,74.75)(56.25,67.5)
\multiput(69.75,34.75)(-.0627907,-.03372093){215}{\line(-1,0){.0627907}}
%\end
%\emline(56.25,67.5)(56.25,67.5)
\put(56.25,27.5){\line(0,1){0}}
%\end
%\emline(70,74.75)(85.5,67.25)
\multiput(70,34.75)(.06950673,-.03363229){223}{\line(1,0){.06950673}}
%\end
\put(73,36.75){\makebox(0,0)[cc]{$m$}}
%\emline(96.5,105.25)(85.5,67.25)
\multiput(96.5,65.25)(-.033639144,-.116207951){327}{\line(0,-1){.116207951}}
%\end
\put(96.75,69.75){\makebox(0,0)[cc]{$v$}}
%\emline(85.5,109.25)(74.5,70.25)
\multiput(85.5,69.25)(-.033639144,-.119266055){327}{\line(0,-1){.119266055}}
%\end
%\emline(74.5,70.25)(85.75,67.25)
\multiput(74.5,30.25)(.1264045,-.0337079){89}{\line(1,0){.1264045}}
%\end
%\emline(56.5,67.5)(74.75,70)
\multiput(56.5,27.5)(.2433333,.0333333){75}{\line(1,0){.2433333}}
%\end
\put(74,27.5){\makebox(0,0)[cc]{$n$}}
\end{picture}

\begin{proof}[Proof of Proposition \ref{muwelldefined}]We argue by induction on $n$. The case $n=1$
follows by Lemma \ref{lem:induct}.

Now let us assume that $n>1$ and that the lemma is true for  partitions of any wall-interval into $n-1$
wall-intervals. Notice first that, according to Lemma~\ref{lem:wallproject}(1) and Lemma \ref{lem:induct}., modulo replacing $(x_i,y_i)$  by its projection with target $(a,b)$, we can assume that the $x_i$'s and $y_i$'s
belong to the interval $I(a,b)$.

We straighten $(a,x_1,y_1,b)$ to  $(a,p_1,q_1,b)$. Then by
Lemma~\ref{lem:wallproject}(2) the sequence $(a,p_1,q_1,b)$ is  geodesic, and we have ${\mathcal W}(x_1\vert
y_1)={\mathcal W}(p_1\vert q_1)$.

By Lemma~\ref{lem:wallgeod} we have ${\mathcal W}(a\vert b)={\mathcal W}(a\vert p_1)\sqcup {\mathcal
W}(p_1\vert q_1)\sqcup {\mathcal W}(q_1\vert b)$. It follows that ${\mathcal W}(a\vert p_1)\sqcup
{\mathcal W}(q_1\vert b)=\sqcup_{i=2}^{n}{\mathcal W}(x_i\vert y_i)$.

We now straighten each path $(a,x_i,y_i,b)$ to  $(a,p_i,q_i,b)$ (when $i>1$). Again we have ${\mathcal W}(x_i\vert y_i)={\mathcal W}(p_i\vert q_i)$ and moreover
$\dist(x_i,y_i)=\dist(p_i,q_i)$ (since $[x_i,p_i,y_i,q_i]$ is a rectangle). Now let us project the  points $p_i$ and $q_i$ onto $I(x,p_1)$ and
$I(q_1,y)$. So set $s_i=m(p_i,x,p_1)$, $t_i=m(q_i,x,p_1)$, $u_i=m(p_i,q_1,y)$ and $v_i=m(q_i,q_1,y)$.

Applying again  Lemma~\ref{lem:wallproject}(1) we get that ${\mathcal W}(p_i\vert q_i)\cap  {\mathcal
W}(a\vert p_1)= {\mathcal W}(s_i\vert t_i)$ and  ${\mathcal W}(p_i\vert q_i)\cap {\mathcal W}(q_1\vert
b)= {\mathcal W}(u_i\vert v_i)$. Thus ${\mathcal W}(p_i\vert q_i) =  {\mathcal W}(s_i\vert t_i)\sqcup
{\mathcal W}(u_i\vert v_i)$, and we get two decompositions: $ {\mathcal W}(a\vert p_1) =
\sqcup_{i=2}^{n}{\mathcal W}(s_i\vert t_i)$ and $ {\mathcal W}(q_1\vert b) = \sqcup_{i=2}^{n}{\mathcal
W}(u_i\vert v_i)$.

If we apply the induction hypothesis to the two decompositions above we see that we are done since
Lemma \ref{aux} ensures that $\dist(p_i,q_i)=\dist(s_i,t_i)+\dist(u_i,v_i)$.
\end{proof}
The following shows that the premeasure satisfies property $(M_1'')$.
\begin{prop}\label{prop:mucontinuous}
Let $(X,\dist)$ be a median space, endowed with convex walls. If $(I_n)_{n\in\naturals}$ is a
non-increasing  sequence of finite disjoint unions of wall-intervals such that $\cap_n I_n=\emptyset$,
then $I_k=\emptyset$ for $k$ large enough.
\end{prop}
\begin{proof}
In what follows we identify a half-space with its
characteristic function. First note that the set of half-spaces bounding a convex wall (i.e. the set of convex subsets whose complement is convex as
well) is a closed subset of $\{0,1\}^X$.  Then the set ${\mathcal H}(x\vert y)$ of half-spaces containing $x$ but not
$y$ is a closed subset of the compact subset of $\{0,1\}^X$ consisting in functions $f:X\to\{0,1\}$
such that $f(x)=1,f(y)=0$. So ${\mathcal H}(x\vert y)$ is compact.

It is enough to argue when $I_0={\mathcal W}(x\vert y)$. Since $(I_n)_{n\in\naturals}$ is non
increasing for each $n$ we have $I_n\subset {\mathcal W}(x\vert y)$. We then define $H_n$ as the set of
half-spaces $h$ such that $\{h,h^c\}\in I_n$, and $x\in h$. It follows that $(H_n)_{n\in\naturals}$ is
non increasing, and has empty intersection.
By projecting onto $I(x,y)$ we have $I_n=\sqcup\ww(x_i\vert y_i)$ for some points $x_i;y_i\in I(x,y)$ (Lemma~\ref{lem:wallproject}(1)).
We know that $\ww(x_i\vert y_i)=\ww(p_i\vert q_i)$ for $p_i=m(x,x_i,y_i),q_i=m(y,x_i,y_i)$, and furthermore $(x,p_i,q_i,y)$ is a geodesic sequence. Thus $H_n=\sqcup\ww(p_i\vert q_i)$ and $H_n$ is compact. It follows that there exists $k$ such that
$H_k=\emptyset$, which implies that $I_k=\emptyset$.
\end{proof}
We now have all the ingredients to finish the proof of Theorem \ref{medwalls}.
\begin{proof}[Proof of Theorem \ref{medwalls}] That the premeasure $\mu$ is well-defined on $\rr$ is the content of
Proposition \ref{muwelldefined}. It obviously satisfies properties $(M_0)$ and  $(M_1')$, while
$(M_1'')$ is proved in Proposition \ref{prop:mucontinuous}.

By Carath\'edory's theorem \ref{thCar}, $\mu^*$ restricted to $\mathcal{A}^*$ is a measure extending
$\mu$, hence its restriction to $\bb$ is also a measure extending $\mu$.

Obviously any isometry of $(X,\dist )$ defines a bijective transformation on $\ww$ preserving $\rr$ and
the premeasure $\mu$, hence the outer measure $\mu^*$ and $\mathcal{A}^*$, hence it defines an
automorphism of the measured space $(\ww , \bb , \mu )$.\end{proof}

%%%%%%%%%%%%%%%
\section{Kernels, median spaces, properties (T) and Haagerup}\label{sect:kernel}
%%%%%%%%%%%%%%%%%%%%%%%%%%%%%%%%%%%%%%%%%%%%%%%%%%%%%%%%%%%%%%%%%%%%%%%%%%%%%%%%%%%%%%

%%%%%%%%%%%%%%%%%%%%%%%%%%%%%%%%%%%%%%%%%%%%%%%%%%%%%%%%%%%%%%
\subsection{Various types of kernels.}
%%%%%%%%%%%%%%%%%%%%%%%%%%%%%%%%%%%%%%%%%%%%%%%%%%%%%%%%%

A \emph{kernel} on a set $X$ is a symmetric map $\psi:X\times X\to \R_+$ such that $\psi(x,x)=0$. For
instance, a pseudo-metric is a kernel.

Let $f:X\to Y$ be a map and let $\phi$ be a kernel on $Y$. The \emph{pull-back of $\phi$ under $f$} is
the kernel $\psi(x,y)=\phi(f(x),f(y))$. Given a class of kernels $\mathcal C$, a kernel $\psi$ on $X$
is \emph{of type $\mathcal C$} if $\psi$ is the pull-back of some kernel in the class $\mathcal C$.

We will be particularly interested in \emph{kernels of median type}, which are obtained by pulling back
a median pseudo-distance. By considering the canonical median metric quotient, we see that any kernel
of median type is also the pull-back of a median distance.

Properties (T) and Haagerup (a-T-me\-na\-bi\-lity) have often been described using conditionally
negative definite kernels, the definition of which we now recall.

\begin{defn}\label{cnd}
A kernel $\psi : X\times X \to \R_+$ is  \emph{conditionally negative definite} if  for every $n\in
\N$, $x_1,...,x_n\in X$ and $\lambda_1,...,\lambda_n\in \R$
    with $\sum_{i=1}^n \lambda_i=0$ the following holds:
    $$\sum_{i=1}^n \sum_{j=1}^n \lambda_i \lambda_j \psi(x_i,x_j)\leq 0\, .$$
\end{defn}
\begin{prop}[\cite{Schoenberg:negdefker}]\label{psialpha}
If $\psi : X\times X \to \R_+$ is a conditionally negative definite kernel and $0<\alpha \leq 1$ then
$\psi^\alpha$ is a conditionally negative definite kernel.
\end{prop}
An example of conditionally negative definite kernel is provided by the following result.
\begin{prop}[\cite{WellsWilliams:Embeddings}, Theorem 4.10]\label{kerlp}
Let $(Y,\bb , \mu)$ be a measured space. Let $0< p\leq 2$, and let $E=L^p (Y,\mu )$ be endowed with the
norm $\|\cdot \|_p$. Then $\psi : E\times E \to \R\, ,\, \psi (x,y)= \|x-y \|_p^p$ is a conditionally
negative definite kernel.
\end{prop}

In some sense, the example in Proposition \ref{kerlp} is universal for conditionally negative definite
kernels, as the following statement shows.

\begin{prop}[\cite{Schoenberg:negdefker}]\label{emb}
A function $\psi : X\times X \to \R_+$ is a conditionally negative definite kernel if and only if there
exists a map $f:X \to H$, where $(H , \|\cdot \| )$ is a Hilbert space, such that
\begin{equation}\label{negdefker}
   \psi (x,y) = \|f(x)-f(y)\|^2\, .
\end{equation}
\end{prop}

The discussion above suggests the following:

\begin{defn}\label{Ddefker}
A function $\psi : X\times X \to \R_+$ is a \emph{kernel of type $p$}, where $0< p\leq 2$, if there
exists a map $f:X \to L^p(Y,\mu)$, for some measured space $(Y,\bb , \mu)$, such that
\begin{equation}\label{pdefker}
   \psi (x,y) = \|f(x)-f(y)\|_p^p\, .
\end{equation}
\end{defn}
\begin{prop}\label{medtype1}A function $\psi : X\times X \to \R_+$ is a kernel of type 1 if and only if it is of median type.\end{prop}

\begin{proof}Since $L^1(Y,\mu)$ is a median space (see Example \ref{exmp:med} (\ref{exmp:L1})), a kernel of type 1 is of median type. Conversely, Corollary \ref{cor:medinL1} shows that a median space embeds in some $L^1(W,\mu)$, so a kernel of median type, by composition with this embedding, will be of type 1.\end{proof}

\begin{rem}\label{rem:pullbCND}
Clearly, the pull-back of a conditionally negative definite kernel (or of a kernel of type $p$) is also
conditionally negative definite (respectively, of type $p$).
\end{rem}
Proposition \ref{emb} states that conditionally negative definite kernels are the same thing as kernels
of type 2. In order to investigate further the relationship between conditionally negative definite
kernels and kernels of type $p$, we recall some results on isometric embeddings of $L^p$--spaces.

%\begin{thm}[\cite{BretDacunhaCastelKriv}, Theorem 1]\label{BDCK1}
%A normed vector space $(E, \|\cdot \|)$ can be embedded linearly and isometrically into $(L^p (X,\mu )
%\, ,\, \|\cdot \|_p)$ for some measured space $(X,\bb , \mu)$ and $1\leq p\leq 2$ if and only if $\|x-y \|^p$ is a conditionally negative definite kernel on $E$.
%\end{thm}

\begin{thm}[Theorems 1 and 7 in \cite{BretDacunhaCastelKriv}]\label{embdlplq}
Let $1\leq p\leq q\leq 2$.
\begin{enumerate}
    \item\label{qinp} The normed space $(L^q (X,\mu ) \, ,\, \|\cdot \|_q)$ can be embedded linearly and isometrically into $(L^p (X',\mu' ) \, ,\, \|\cdot
\|_p)$ for some measured space $(X',\bb' , \mu')$.
    \item\label{pinq} If $L^p (X,\mu )$ has infinite dimension then $(L^p (X,\mu ) \, ,\, \|\cdot \|_p^\alpha)$ can
    be embedded isometrically into $(L^q (X',\mu' ) \, ,\, \|\cdot
\|_q)$ for some measured space $(X',\bb' , \mu')$ if and only if $0<\alpha \leq \frac{p}{q}$.
\end{enumerate}
\end{thm}

\begin{rem}
Note that according to \cite{JohnsonRandrian}, the space $l_p$
with $p>2$ does not coarsely embed into a Hilbert space.
\end{rem}

\begin{rem}\label{rem:typepsubm}
Theorem \ref{embdlplq}, (\ref{qinp}), implies that every metric space that can be isometrically
embedded in a space $L^p (X,\mu )$ with $p\in [1,2]$ (\emph{metric space of type p} in the terminology
of \cite{BretDacunhaCastelKriv} and \cite{FarautHarz}) is a submedian space. See \cite{FarautHarz} for
examples of such spaces.
\end{rem}

Using these results we can now establish a more precise relationship between kernels of type $p$ and
conditionally negative definite.

\begin{cor}\label{cndtypep}
\begin{enumerate}
    \item If $\psi$ is a kernel of type $p$ for some $0< p\leq 2$ then $\psi$
    is a conditionally negative definite kernel.
    \item If $\psi$ is a conditionally negative definite kernel and $1\leq p\leq 2$ then
    $\psi^{\frac{p}{2}}$ is a kernel of type $p$.
\end{enumerate}
\end{cor}

\begin{proof}Let $X$ be an arbitrary space and let $\psi :X \times X \to \R_+$.

\medskip

(1) follows from Proposition \ref{kerlp} and Remark \ref{rem:pullbCND}.

 \medskip

(2) According to Proposition \ref{emb} there exists a map $g : X \to L^2(X,\mu )$ such that $\psi
(x,y)= \|g(x)-g(y)\|_2^2$. By Theorem \ref{embdlplq}, (1), there exists an isometric embedding $F:
\left( L^2(X,\mu)\, ,\, \|\cdot \|_2 ), \right)\to (L^p (X',\mu' ) \, ,\, \|\cdot \|_p)$. Consequently
$\psi (x,y)= \|g(x)-g(y)\|_2^2 = \|F(g(x))-F(g(y))\|_p^2\, $, and $\psi^{p/2}$ is a kernel of type $p$.
\end{proof}

\begin{rem}\label{rem:subm}
\begin{enumerate}
    \item\label{spaces} By Proposition \ref{medtype1}, Corollary \ref{cndtypep} and Proposition \ref{emb},
     every submedian space $(X,\dist )$ has the
property that $(X,\dist^{1/2})$ can be embedded isometrically in a Hilbert space. This can be refined
(\cite{Assouad:L1}, \cite[Proposition 2.5]{DezaGrishLaurent:survey}) to the sequence of implications:
$(X, \dist)$ submedian $\Rightarrow$ $(X, \dist)$ hypermetric $\Rightarrow$ $(X,\dist^{1/2})$
spherically $L^2$-embeddable $\Rightarrow$ $(X,\dist^{1/2})$ $L^2$-embeddable.

Recall that a kernel $\psi : X\times X \to \R$ (in particular a metric) is \emph{hypermetric} if for
any finite sequence $x_1,...,x_n$ in $X$ and any integers $\lambda_1,...\lambda_n$ such that
$\sum_{i=1}^n \lambda_i=1$, we have $\sum_{i,j=1}^n \lambda_i \lambda_j \psi (x_i,x_j)\leq 0$. A kernel
is \emph{spherical} if its restriction to any finite subset of $X$ coincides with a pull-back of a
metric on an Euclidean unit sphere. A metric space is called \emph{spherically $L^2$-embeddable} if its
distance is a spherical kernel.

It follows that any submedian space  $(X,\dist )$ has the property that all its finite subsets endowed
with the metric $\dist^{1/2}$ are isometric to subsets on an Euclidean unit sphere. This holds even for
submedian spaces of negative curvature, like $\hyperbolic^n_\R$, as was first noticed by Robertson in
\cite[Corollary 3.2]{Robertson:Crofton}.
    \item\label{kernels} The above implications can be reformulated in terms of kernels thus:
  $\psi$  kernel of type $1$ $\Rightarrow$ $\psi$ hypermetric kernel $\Rightarrow$ $\psi^{1/2}$
  spherical kernel $\Rightarrow$
   $\psi$ kernel of type $2$.
\end{enumerate}

\end{rem}
Robertson and Steger defined in \cite{RobertsonSteger:negdefker} an alternate type of kernels.
\begin{defn}[Robertson and Steger \cite{RobertsonSteger:negdefker}]A \emph{measure definite kernel}
on a space $X$ is a map $\psi : X\times X \to \R_+$ such that there exists a measured space $(\mm ,
\bb, \mu)$ and a map $S:X\to \bb$, $x\mapsto S_x$, satisfying $\psi (x,y)= \mu (S_x \vartriangle
S_y)$.\end{defn}
In \cite{RobertsonSteger:negdefker} it is asked (Question (i)) whether measure definite kernels can be
given an intrinsic characterization among the conditionally negative definite kernels. It turns out
that measure definite kernels are very much related to structures of space with measured walls, as well
as to median spaces (see Lemma \ref{kerm}). This relationship allows us to answer this question
(Corollary \ref{mesdefintr}).
\begin{lem}\label{kerm}
A kernel $\psi : X\times X \to \R_+$ on a space $X$ is  measure definite  if and only if it is of median type,
in other words it is the pull-back of a median metric.
{Moreover when $X$ is a topological space the kernel $\psi$ is continuous if and only if the pull-back map
$f$ is continuous.}
\end{lem}

\begin{proof} Assume that $\psi$ is a \mdk on $X$. Then there exists a map $S:X\to \bb$, $x\mapsto S_x$, where
$(\mm , \bb, \mu)$ is a measured space, and $\psi (x,y)= \mu (S_x \vartriangle S_y)$. Fix some base point $x_0$ and  endow $\bb_{S_{x_0}}$ with the structure of median pseudo-metric
space described in Example~\ref{exmp:med}(\ref{exmp:symdiff}). Then $\psi$ is the pull-back under $S$
of this median pseudo-metric.

Conversely, consider a map $f$ from $X$ to a median space $(Y,\dist)$ such that $\psi (x,x')= \dist
(f(x),f(x'))$. By Theorem \ref{medwalls}, there exists a set of convex walls $\ww$ on $Y$, a
$\sigma$-algebra ${\mathcal B}$ on $\ww$ and a measure $\mu$ on ${\mathcal B}$ such that the 4-tuple
$(Y,\ww,{\mathcal B},\mu)$ is a space with measured walls, and moreover $\dist (y,y') = \mu (\ww
(y|y'))$.

We fix a point $x_0$ in $X$ and we define the map $S :X \to \bb\, ,\, S_x= \ww (f(x) | f(x_0))$. Then
$\mu (S_a \vartriangle S_b)= \mu \left( \ww (f(a) | f(x_0)) \vartriangle \ww (f(b) | f(x_0)) \right)=
\mu (\ww (f(a)|f(b)))=\dist (f(a),f(b))=\psi (a,b)$.

Obviously $f$ continuous implies $\psi$ continuous. Conversely, assume that $\psi$ is continuous. If
$y\in X$ is close to $x\in X$ then $(x,y)$ is close to $(x,x)$ hence $\psi (x,y)=\dist (f(x),f(y))$ is
close to $\psi (x,x)=0$.
\end{proof}

The following statement is an improvement of \cite[Proposition 1.2]{RobertsonSteger:negdefker} and of
\cite[Proposition 2]{CherixMartinValette}.

\begin{lem}\label{kermeas}
A map $\psi : X\times X \to \R_+$ on a space $X$ is a measure definite kernel if and only if there
exists a structure of space with measured walls $(X, \ww , \bb , \mu )$ on $X$ such that $\psi (x,x')=
\mu (\ww(x|x'))$.
\end{lem}

\proof The if part follows from Lemmata \ref{lem:embedL1} and \ref{kerm}.

Conversely, assume that $\psi$ is a measure definite kernel on $X$. By Proposition \ref{kerm} the
kernel $\psi$ is the pull-back of a median distance: $\psi(x,y)=\dist(f(x),f(y))$ for some map $f:X\to
Y$ where $(Y,\dist)$ is a median space.  Consider the structure of space with measured walls on $Y$
given by Theorem \ref{medwalls}. The pull-back structure of space with measured walls on $X$ has $\psi$
as wall pseudo-distance, according to Lemma~\ref{lem:pullback}.
\endproof

\begin{prop}\label{mesdefconddefI}
A kernel is measure definite if and only if it is of type 1.
\end{prop}

\proof Follows directly from Lemma~\ref{kerm} and Proposition \ref{medtype1}.
\endproof
\begin{cor}\label{mesdefintr}
A kernel $\psi : X\times X \to \R$ is measure definite if and only if  $\psi$ satisfies the triangular
inequality, moreover for every finite subset $F$ in $X$, $\psi|_{F\times F}$ is equal to
    $\sum_{S\subseteq F}\lambda_S \delta_S$ for some $\lambda_S\geq 0$, where $\delta_S(x,y)=1$ if
    $\psi(x,y)>0$ and $S \cap \{x,y\}$ is of cardinality $1$, $\delta_S(x,y)=0$ otherwise.
\end{cor}

\proof It follows immediately from Proposition \ref{mesdefconddefI}, from the fact that a metric space
is isometrically embeddable into an $L^1$-space if any finite subset of it is
(\cite{Assouad:analytique}, \cite{AssouadDeza}), and from Remark \ref{finite}.\endproof
\begin{cor}\label{mesdefconddef}
\begin{enumerate}
    \item  Every measure definite kernel is conditionally negative definite.

    \me

    \item If $\psi$ is a conditionally negative definite kernel then $\sqrt{\psi
    }$ is a measure definite kernel.
\end{enumerate}
\end{cor}

Statement (1) in Corollary \ref{mesdefconddef} has already been proved in
\cite{RobertsonSteger:negdefker}, where it appears as Proposition 1.1, while statement (2) has been
proved in \cite[Proposition 1.4(i)]{RobertsonSteger:negdefker} under the extra assumption that the set on
which the kernel is defined is countable.

%%%%%%%%%%%%%%%%%%%%%%%%%%%%%%%%%%%%%%%%%%%%%%%%%%%%%%%%%%%%%%
\subsection{Properties (T) and Haagerup and actions on median, measured walls and $L^p$--spaces.}\label{kerThag}
%%%%%%%%%%%%%%%%%%%%%%%%%%%%%%%%%%%%%%%%%%%%%%%%%%%%%%%%%

\begin{defn}
A function $\Phi :G \to \R_+$ defined on a group is\emph{ conditionally negative definite} if the
function $G\times G \to \R_+\, ,\, (g,h)\mapsto \Phi (g^{-1}h)$, is a conditionally negative definite
kernel.
\end{defn}
Recall that a function $\Phi$ is called \emph{proper} if $\lim_{g\to \infty } \Phi (g)=\infty$. Here
$g\to \infty$ means that $g$ leaves any compact subset.

If a conditionally negative definite kernel $\psi : G\times G \to \R_+$ is \emph{left invariant}, i.e.
$\psi (g_1,g_2)= \psi (hg_1,hg_2)$ for every $h,g_1,g_2$ in $G$, then the map $\Phi :G \to \R_+$
defined by $\Phi (g)= \psi (1,g)$ is a conditionally negative definite function. If $\Phi$ is proper we
say that the kernel $\psi$ is \emph{proper}.

We also recall that a \emph{second countable space} is a topological space satisfying the second axiom
of countability, that is such that its topology has a countable base. A second countable space is
\emph{separable} (i.e. has a countable dense subset) and \emph{Lindel\"of} (i.e. every open cover has a
countable sub-cover). The converse implications do not hold in general, but they do for metric spaces.

Characterizations of properties (T) and Haagerup (also called a-T-me\-na\-bi\-lity) using conditionally
negative definite kernels are well-known
 and can be found in the literature. We recall here the relevant ones.
\begin{thm}[\cite{Delorme}, \cite{Guichardet}, \cite{AkermannWalter}, \cite{delaHarpeValette:proprieteT},
\cite{CCJJV}]\label{THker} Let $G$ be a second countable, locally compact group.
\begin{enumerate}
    \item  The group $G$ has property (T) if and only if every
    continuous conditionally negative definite function on $G$ is bounded (equivalently, every
    continuous left invariant conditionally negative definite kernel on $G$ is bounded).
    \item   The group $G$ has the Haagerup property if and only if there exists a
    continuous proper conditionally negative definite function on $G$ (equivalently, there exists a
    continuous proper left invariant conditionally negative definite kernel on $G$).
\end{enumerate}
\end{thm}

Theorem \ref{THker} and Corollary \ref{cndtypep} imply the following.

\begin{cor}\label{c1}
Let $G$ be a second countable, locally compact group.
\begin{enumerate}
    \item  If the group $G$ has property (T) then for every $p\in (0,2]$, every
    continuous left invariant kernel of type $p$ on $G$ is bounded.
    \item   The group $G$ has the Haagerup property if for some $p\in (0,2]$,
     there exists a continuous proper left invariant kernel of type $p$ on $G$.
\end{enumerate}
\end{cor}
\begin{rem}\label{rem:c1}
For $p\in [1,2]$ the converse statements in Corollary \ref{c1} immediately follow from Corollary
\ref{cndtypep}, (2).
\end{rem}
Corollary \ref{c1} can be reformulated in terms of actions of the group on subsets of $L^p$-spaces, as
follows.
\begin{cor}[\cite{Delorme}, \cite{AkermannWalter}, \cite{WellsWilliams:Embeddings}]\label{c2}
Let $G$ be a second countable, locally compact group.
\begin{enumerate}
    \item\label{t}  If $G$ has property (T) then for every $p\in (0,2]$, every
    continuous action by isometries of $G$ on a subset of a space $L^p(X,\mu)$ has bounded orbits.
    \item\label{at}   The group $G$ has the Haagerup property if there exists  $p\in (0,2]$,
     and a continuous proper action by isometries of $G$ on a subset of some $L^p(X,\mu)$.
\end{enumerate}
\end{cor}

\proof Both (\ref{t}) and (\ref{at}) follow from the fact that if $S$ is a subset of some $L^p(X,\mu)$
and there exists an action of $G$ on $S$ by isometries, $G\times S \to S\, ,\, (g,s)\mapsto g\cdot s$,
then for any $s\in S$ the map $\psi (g,h)= \| g\cdot s - h\cdot s\|_p^p$ is a continuous left invariant
kernel of type $p$ on $G$.
\endproof

\begin{rem}\label{nowak}
In \cite{PNowak:lp} the following result is stated: a second countable
locally compact group has the Haagerup property if and only if for some
(for all) $p\in (1,2)$ the group has a proper affine isometric action on
$L^p[0,1]$. The proof in that paper has been completed in an updated version of his preprint appearing on arXiv at
\cite{PNowak:arxiv}.
\end{rem}

\me

The converse statements in Corollary \ref{c2} (and their stronger versions, with ``every $p\in (0,2]$''
replaced by ``there exists $p\in (0,2]$'' in (\ref{t}), and the opposite replacement done in
(\ref{at})) follow immediately from the following fact. Given $\psi$ a continuous proper left invariant
kernel of type $p$ on $G$, that is a map $\psi :G\times G\to \R_+$ defined by $\psi (g,h)=\| f(g)
-f(h)\|_p^p$, where $f:G\to L^p(X,\mu)$ is continuous, one can define a continuous action by isometries
of $G$ on $f(G)$ by $g\cdot f(h)=f(gh)$.

Much stronger versions of the converse statements in Corollary \ref{c2} are provided by Corollary
\ref{c3}. But in order to obtain those, we first need to obtain improved converse statements for $p=1$.
Indeed, for this value of $p$, the sufficient condition to have property (T) can be weakened: it
suffices to look at actions of $G$ on median subspaces of $L^1$--spaces. Also, Haagerup property
implies more for $p=1$: the existence of a continuous proper action by isometries of $G$ on a median
subspace of some $L^1$--space. Both statements are straightforward consequences of the following
result.
\begin{thm}\label{kermgrI}
Let $G$ be a separable topological group.
\begin{enumerate}
  \item\label{thaction}   {If $G$ acts continuously by isometries on a
median space $(X,\dist)$ and $x\in X$ then $\psi : G\times G \to \R_+$, $\psi (g,g')= \dist (g\cdot x,g'\cdot x)$ is a continuous left
invariant kernel of type 1.}
  \item\label{tker}   {If $\psi : G\times G \to \R_+$ is a continuous left
invariant kernel which is the square root of a kernel of type 2 (hence $\psi$ is a kernel of type 1) then there exists a continuous action by isometries of $G$ on a
median space $(X,\dist)$, and a point $x\in X$ such that}
   {$$\psi (g,g')= \dist (g\cdot x,g'\cdot x)\, .$$}
\end{enumerate}
\end{thm}

We first need to establish equivariant versions of Lemmata \ref{kerm} and \ref{kermeas} when $X$ is a
group $G$. In the particular case when the group  is countable, Lemma \ref{kermeas} has the following
equivariant version.
\begin{prop}[\cite{RobertsonSteger:negdefker}]\label{kermeasgr}
  {Let $\Gamma$ be a countable group.}
\begin{enumerate}
  \item\label{measure}   {If $\Gamma$ is endowed with a left invariant structure of space with measured walls $(\Gamma, \ww , \bb , \mu )$  then  $\psi : \Gamma\times \Gamma \to \R_+$ defined by $\psi (g,g')= \mu (\ww(g|g'))$ is a left invariant measure definite
kernel.}
  \item\label{kernel}   {If $\psi : \Gamma\times \Gamma \to \R_+$ is the square root of a
left invariant conditionally negative definite kernel, then $\Gamma$ can be endowed with a left invariant structure of space with measured walls $(\Gamma, \ww , \bb , \mu )$ such that $\psi (g,g')= \mu (\ww(g|g'))$.}
\end{enumerate}
\end{prop}
\begin{proof}   {(\ref{measure}) follows immediately from the definition of a measure definite
kernel. It appears in \cite{RobertsonSteger:negdefker} as Proposition 1.1.}

  {(\ref{kernel}) follows from Proposition 1.4 and the
proof of Theorem 2.1 on p.252 in \cite{RobertsonSteger:negdefker}.}\end{proof}

This implies the following equivariant version of Lemma \ref{kerm}.
\begin{lem}\label{kermgr}
Let $\Gamma$ be a countable group.
\begin{enumerate}
  \item\label{action}   {If $\Gamma$ acts by isometries on a median
space $(X,\dist)$, and $x$ is a point in $X$ then $\psi : \Gamma\times \Gamma \to \R_+$,  $\psi (g,g')= \dist (g\cdot x,g'\cdot x)$ is a left invariant kernel of type $1$.}
  \item\label{kernmed}   {If $\psi : \Gamma\times \Gamma \to \R_+$ is a left invariant kernel which is the square root of a conditionally negative definite kernel (hence $\psi$ is of type $1$) then there exists an action by isometries of $\Gamma$ on a median
space $(X,\dist)$, and a point $x\in X$ such that $\psi (g,g')= \dist (g\cdot x,g'\cdot x)$.}
\end{enumerate}
\end{lem}
\begin{proof}[Proof of Theorem \ref{kermgrI}]
(\ref{thaction}) follows from the fact that median spaces isometrically embed in $L^1$-spaces by
Corollary~\ref{cor:medinL1}.

(\ref{tker}) Let $\Gamma$ be a countable dense subgroup in $G$.   {Restrict $\psi$ to
$\Gamma$ and apply Lemma \ref{kermgr}, (\ref{kernmed}):} there exists an action by
isometries of $\Gamma$ on a median space $(X, \dist )$ and $x\in X$ such that $\psi (\gamma ,\gamma')=
\dist (\gamma \cdot x,\gamma'\cdot x)$. The metric completion of $(X,\dist)$  is still median by
Proposition \ref{cor:medcomplete}, and any isometry of $X$ extends uniquely to an isometry of the
completion. We get an action of $\Gamma$ on a complete median space that still induces the kernel
$\psi$. Thus we may - and will - assume that the median space $(X,\dist)$ is already complete.

 The map $f:\Gamma \to X $ sending $
\gamma$ to $ \gamma \cdot x$ is uniformly continuous since $\psi$ is continuous (we endow $\Gamma\subset G$ with
the induced topology). Since $X$ is complete it follows that $f$ extends to a continuous map $G\to X$
still verifying $\psi (g,g')= \dist (f(g),f(g'))$ (for all $g,g'$ in $G$). As usual the left
invariance of $\psi$ implies that $g\cdot f(h)=f(gh)$ defines an action of $G$ by isometries on $f(G)$.
And the continuity of $\psi$ implies that the action is continuous.

To end the argument it suffices to prove that the above action of $G$ on $f(G)$ extends to an action by
isometries on the median hull of $f(G)$ in $X$, which we denote by $M$. Set $M_0=f(G)$ and then define
inductively $M_{i+1}=\{m(x,y,z)\mid (x,y,z)\in (M_i)^3\}$. Clearly $M$ is the ascending union of the
$M_i$'s.   {Every element $g\in G$ defines an isometry $g:M_0\to M_0$. We first note that there is at most one isometric embedding $\overline g:M\to X$ extending $g:M_0\to M_0$, and $\overline g(M)\subset M$. Indeed, since isometries commute with the median
map,  $\overline g$ is completely determined on $M_1$ and $\overline g(M_1)\subset M_1$, then
$\overline g$ is completely determined on $M_2$ and $\overline g(M_2)\subset M_2$, and so on.}  We now
prove the existence of such an isometric extension.

 Choose a sequence $\gamma_n$ of elements of $\Gamma$ converging to $g$ in the topological group $G$. Then for every $f(h)\in M_0$, $\gamma_n f(h)= f(\gamma_n h)$ converges to $f(gh)=gf(h)$ by continuity of $f$. We prove by induction on $i$ that for any $m\in M_i$ the sequence $\gamma_n(m)$ converges. This is true for $m\in M_0$. Assume we know that $\gamma_n(p)$ converges for every $p\in M_i$, and let $m$ denote an element of $ M_{i+1}$. Write $m=m(x,y,z)$ with $(x,y,z)\in (M_i)^3$.
  Since $\gamma_n$ acts by isometry on the whole space
$X$ we have $\gamma_n(m)=m(\gamma_n(x),\gamma_n(y),\gamma_n(z))$. By induction the three sequences $(\gamma_n(x))_{n\ge 0},(\gamma_n(y))_{n\ge 0}(\gamma_n(z))_{n\ge 0}$ are convergent.
By the continuity of the median map (Corollary~\ref{cor:medlip}) it follows that $(\gamma_n(m))_{n\ge 0}$ converges.
Denote $\overline g:M\to X$ the pointwise limit of $\gamma_n$ on $M$.
It immediately follows that $\overline g$ is an isometric embedding which extends $g:M_0\to M_0$.
By the remarks above we have that $\overline g(M)\subset M$.

Using the uniqueness of the extension it is now straightforward to check that the maps $ \overline g$
are isometries of $M$ (with inverse  $\overline{g^{-1}}$), and finally that $g\mapsto \overline g$
defines an action of $G$ by isometries on $M$ that extends the action of $G$ on $f(G)$, and thus still
induces the kernel $\psi$.\end{proof}

Theorem \ref{kermgrI} allows to obtain some results concerning a structure of space with measured walls
on the complex hyperbolic space.

\begin{cor}[walls in the complex hyperbolic space]\label{cor:complex}
The complex hyperbolic space $\field H^n_\C$ admits a structure of space with measured walls such that:
\begin{enumerate}
    \item the induced wall metric is $\dist^{1/2}$, where $\dist$ is the hyperbolic distance;
    \item the walls are all the convex walls with respect to the metric $\dist^{1/2}$;
    \item $SU(n,1)$ acts by isomorphisms on this structure.
\end{enumerate}

\end{cor}

\proof According to \cite{FarautHarz} the complex hyperbolic space $\field H^n_\C$ equipped with the
metric $\dist^{1/2}$ can be embedded into a Hilbert space. It follows, by Theorem \ref{embdlplq},
(\ref{qinp}), that $\psi : \field H^n_\C \times \field H^n_\C \to \R_+ \, ,\, \psi (x,y)=\dist^{1/2}
(x,y)\, ,$ is   {the square root of a kernel of type $2$}, in the terminology of Definition \ref{Ddefker}. Obviously $\psi$ is continuous and left-invariant with respect to the action of  $G=SU(n,1)$.

 Via the identification of $\field
H^n_\C$ with $G/K$, where $K=SU(n)$, the kernel $\psi$ induces a left invariant pull-back kernel
$\psi_G : G\times G \to \R_+$. Theorem \ref{kermgrI} implies that $G$ acts by isometries on a median
space $(X, \dist_X)$ such that $\psi_G(g,g')=\dist_X(g\cdot x, g'\cdot x)$ for some $x\in X$. It
follows easily that the map $g\mapsto g\cdot x$ factorizes to a $G$-equivariant isometric embedding $gK
\mapsto gx$ of $(\field H^n_\C , \dist^{1/2})$ into $(X, \dist_X)$. All the required statements then
follow from Lemma \ref{lem:pullback} and from Theorem \ref{tsepp}.
\endproof

We now prove the results stated in the introduction.
\begin{proof}[Proof of Theorem \ref{impl1}] By Corollary~\ref{c1} and Remark \ref{rem:c1},
 property (T) and a-T-me\-na\-bi\-lity for a group $G$ are characterized by
 properties of continuous left invariant kernels of type 1.
   {By Theorem \ref{kermgrI}, (\ref{thaction}), continuous actions of
 $G$ on median spaces induce such kernels. On the other hand, a kernel of type 1 is by Corollary \ref{cndtypep} also of type 2, hence its square root is defined by a continuous action on a median space, according to Theorem \ref{kermgrI}, (\ref{tker}). Theorem \ref{impl1} follows,} since bounded kernels correspond
 to actions with bounded orbits, and proper kernels correspond to proper actions.
  \end{proof}

\begin{proof}[Proof of Theorem \ref{implmw}]
If a group acts continuously on a median space by isometries then the group acts continuously by
automorphisms on the structure of measured walls associated to it, by Theorem \ref{medwalls}. This and
Theorem \ref{impl1} give the direct implication in (2) and the converse implication in (1).

 {On the other hand, a space with measured walls is a submedian space, by Corollary \ref{cor:premedian}, hence a subspace of an $L^1$--space by Corollary \ref{cor:medinL1}. Corollary \ref{c2} then gives the direct implication in (1) and the converse implication in (2).}
\end{proof}

\begin{proof}[Proof of Corollary \ref{c3}] A continuous action of a group $G$
on a space with measured walls $(X, \ww , \mu )$ induces by Lemma \ref{lem:actlp} a continuous action
by affine isometries on $L^p (\hh , \mu_\hh )$ for any $p>0$, defined by $g\cdot f = \pi_p(g) (f) +
\chi_{\sigma_{gx}} - \chi_{\sigma_x}$, where $x$ is an arbitrary point in $X$.

The hypothesis in case (1) implies that the orbit of the constant function zero, composed of the
functions $\chi_{\sigma_{gx}} - \chi_{\sigma_x},\, g\in G$, is bounded. This implies that the orbit of
$x$ is bounded. It remains to apply Theorem \ref{implmw}, (1).

(2) If $G$ is a-T-menable then we may assume by  Theorem \ref{implmw}, (2), that $G$ acts on $(X, \ww ,
\bb , \mu )$ such that $\pdist_\mu(x,gx)\to \infty$ when $g\to \infty $. Hence the action of $G$ on
$L^p (\hh , \mu_\hh )$ is proper.\end{proof}

\me

\begin{proof}[Proof of Theorem \ref{Tdyn}] The only if part of (1) and the if part of (2) immediately
follow from Corollary \ref{c1} and Proposition~\ref{mesdefconddefI}.

The if part of (1) and the only if part of (2) follow from Theorem \ref{impl1} and from Corollary
\ref{cor:medinL1}.
\end{proof}

%\begin{quest}
%Can a theorem similar to Theorem \ref{impl1} be formulated for subspaces with special properties of an
%$L^p$-space, $1<p<2$ ?

%For instance, what are the metric properties of the images of
%median subspaces of an $L^1$-space under the embedding described
%in Theorem \ref{embdlplq}, (2) ?
%\end{quest}

%%%%%%%%%%%%%%%%%%%%%%%%%%%%%%%%
%\bibliographystyle{amsalpha} %
%\bibliography{cornelia} %

\begin{thebibliography}{BFGM07}


\bibitem[AW81]{AkermannWalter}
C.A. Akemann and M.E. Walter, \emph{{Unbounded negative definite functions}},
  Canad. J. Math. \textbf{33} (1981), 862--871.

\bibitem[Alp82]{Alperin:propT}
R.~Alperin, \emph{{Locally compact groups acting on trees and Property T}},
  Monatsh. Math. \textbf{93} (1982), 261--265.

\bibitem[AD82]{AssouadDeza}
P.~Assouad and M.~Deza, \emph{{Metric subspaces of $L^1$}}, Publications
  math\'ematiques d'Orsay \textbf{3} (1982).


\bibitem[Ass]{Assouad:analytique}
Patrice Assouad, \emph{Plongements isom\'etriques dans {$L\sp{1}$}: aspect
  analytique}, Initiation {S}eminar on {A}nalysis: {G}. {C}hoquet-{M}.
  {R}ogalski-{J}. {S}aint-{R}aymond, 19th {Y}ear: 1979/1980, Publ. Math. Univ.
  Pierre et Marie Curie, vol.~41, Univ. Paris VI, Paris, pp.~Exp. No. 14, 23.

\bibitem[Ass81]{Assouad:analytiqueComment}
\bysame, \emph{Comment on: ``{I}sometric embeddings in {$L\sp{1}$}: analytic
  aspect''}, Initiation {S}eminar on {A}nalysis: {G}. {C}hoquet-{M}.
  {R}ogalski-{J}. {S}aint-{R}aymond, 20th {Y}ear: 1980/1981, Publ. Math. Univ.
  Pierre et Marie Curie, vol.~46, Univ. Paris VI, Paris, 1981, pp.~Comm. No.
  C4, 10.

\bibitem[Ass84]{Assouad:L1}
P.~Assouad, \emph{{Sur les in\'egalit\'es valides dans $L^1$}}, Europ. J.
  Combinatorics \textbf{5} (1984), 99--112.


\bibitem[BFGM07]{BFGM}
U.~Bader, A.~Furman, T.~Gelander, and N.~Monod, \emph{{Property (T) and
  rigidity for actions on Banach spaces}}, Acta Math. \textbf{198} (2007),
  57--105.

\bibitem[BS97]{BallmannSwiat:polyhedra}
W.~Ballmann and J.~\'Swi\c{a}tkowski, \emph{{On $L^2$-cohomology and property
  (T) for automorphism groups of polyhedral cell complexes}}, Geom. Funct.
  Anal. \textbf{7} (1997), 615--645.

\bibitem[BH83]{BandeltHedlikova}
H.-J. Bandelt and J.~Hedlikova, \emph{Median algebras}, Discrete Math.
  \textbf{45} (1983), 1--30.

\bibitem[BC08]{ChepoiBandelt}
H.-J.~Bandelt and V.~Chepoi, \emph{Metric graph theory and
  geometry: a survey}, Surveys on discrete and computational geometry, Contemp.
  Math., vol. 453, Amer. Math. Soc., Providence, RI, 2008, pp.~49--86.



\bibitem[Bas01]{Basarab}
S.~A. Basarab, \emph{The dual of the category of generalized trees}, {Ann.
  \c{S}tiin\c{t}. Univ. Ovidius Constan\c{t}a Ser. Mat.} \textbf{9} (2001),
  1--20.


\bibitem[Bau01]{HBauer:measure}
H.~Bauer, \emph{{Measure and Integration Theory}}, Studies in Math., vol.~26,
  de Gruyter, Berlin, 2001.


\bibitem[BdlHV]{BekkaHarpeValette}
B.~Bekka, P.~de~la Harpe, and A.~Valette, \emph{{Kazhdan's Property
  (T)}}, New Mathematical Monographs \textbf{11}, Cambridge University Press, Cambridge, 2008.

\bibitem[BdlHV01]{Oberwolfach}
B.~Bekka, P.~de~la Harpe, and A.~Valette, \emph{{Mini-Workshop: Geometrization
  of Kazhdan's Property (T) (July 2001)}}, Report no. 29/2001, Mathematisches
  Forschungsinstitut Oberwolfach, 2001.

\bibitem[BF95]{BestvinaFeighn:stableactions}
M.~Bestvina and M.~Feighn, \emph{Stable actions of groups on real trees},
  Invent. Math. \textbf{121} (1995), no.~2, 287--321.


\bibitem[BK47]{BirkhoffKiss}
G.~Birkhoff and S.~A. Kiss, \emph{A ternary operation in distributive
  lattices}, Bull. Amer. Math. Soc. \textbf{53} (1947), 749--752.


\bibitem[BP03]{BourdonPajot:cohomologie}
M.~Bourdon and H.~Pajot, \emph{{Cohomologie $l_p$ et espaces de Besov}}, J.
  Reine Angew. Math. \textbf{558} (2003), 85--108.

\bibitem[BJS88]{BozejkoJS}
M.~Bo\.{z}ejko, T.~Januszkiewicz, and R.~Spatzier, \emph{{Infinite Coxeter
  groups do not have Property (T)}}, J. Op. Theory \textbf{19} (1988), 63–--67.


\bibitem[Bo{\.z}89]{Bozejko:trees}
M.~Bo{\.z}ejko, \emph{Positive-definite kernels, length functions on groups
  and a noncommutative von {N}eumann inequality}, Studia Math. \textbf{95}
  (1989), no.~2, 107--118.


\bibitem[BCK66]{BretDacunhaCastelKriv}
J.~Bretagnolle, D.~Dacunha Castelle, and J.-L. Krivine, \emph{{Lois stables et
  espaces $L^p$}}, Ann. IHP, section B \textbf{2} (1966), no.~3, 231--259.


\bibitem[CN05]{ChatterjiNiblo}
I.~Chatterji and G.~Niblo, \emph{From wall spaces to {$\rm CAT(0)$}
  cube complexes}, Internat. J. Algebra Comput. \textbf{15} (2005), no.~5-6,
  875--885.

  \bibitem[CDH07]{CDH-geometry}
I.~Chatterji, C.~Dru\c{t}u, and F.~Haglund, \emph{{Geometry of median spaces}},
  preprint, 2007.

\bibitem[Che00]{Chepoi:graphs}
V.~Chepoi, \emph{{Graphs of Some CAT(0) Complexes}}, Adv. in Appl. Math.
  \textbf{24} (2000), 125--179.

\bibitem[CCJ{\etalchar{+}}01]{CCJJV}
P.A. Cherix, M.~Cowling, P.~Jolissaint, P.~Julg, and A.~Valette, \emph{{Groups
  with the Haagerup Property (Gromov's a-T-menability)}}, Progress in Math.,
  vol. 197, Birkhauser, Boston, 2001.


\bibitem[CMV04]{CherixMartinValette}
P.A. Cherix, F.~Martin, and A.~Valette, \emph{Spaces with measured walls, the
  {H}aagerup property and property {(T)}}, Ergod. Th. \& Dynam. Sys.
  \textbf{24} (2004), 1895--1908.


\bibitem[dCSV07]{CornStaldValette:lamp}
Y.~de~Cornulier, Y.~Stalder, and A.~Valette, \emph{{Proper actions of
  lamplighter groups associated with free groups}}, preprint
  {\textsc{arXiv:math.GR/07072039}}, 2007.

\bibitem[dCTV06]{CornTessValette:lp}
Y.~de~Cornulier, R.~Tessera, and A.~Valette, \emph{{Isometric group actions on
  Banach spaces and representations vanishing at infinity}}, preprint
  \textsc{arXiv:math.RT/0612398}, 2006.


\bibitem[dlHV89]{delaHarpeValette:proprieteT}
P.~de~la Harpe and A.~Valette, \emph{{La propri\'et\'e (T) de Kazhdan pour les
  groupes localement compacts}}, Ast\'erisque, vol. 175, Soci\'et\'e
  {M}ath\'ematique de {F}rance, 1989.


\bibitem[Del77]{Delorme}
P.~Delorme, \emph{{$1$-cohomologie des repr\'esentations unitaires des groupes
  de Lie semi-simples et r\'esolubles. Produits tensoriels et
  repr\'esentations}}, Bull. Soc. Math. France \textbf{105} (1977), 281--336.

\bibitem[DGL95]{DezaGrishLaurent:survey}
M.~Deza, V.~P. Grishukhin, and M.~Laurent, \emph{Hypermetrics in the geometry
  of numbers}, DIMACS Series in Discrete Mathematics and Computer Science
  \textbf{20} (1995).


\bibitem[FH74]{FarautHarz}
J.~Faraut and K.~Harzallah, \emph{Distances hilbertiennes invariantes sur un
  espace homog\`ene}, Ann. Institut Fourier \textbf{3} (1974), no.~24,
  171--217.

\bibitem[Ger97]{Gerasimov:semisplittings}
V.~N. Gerasimov, \emph{Semi-splittings of groups and actions on cubings},
  Algebra, geometry, analysis and mathematical physics (Russian) (Novosibirsk,
  1996), Izdat. Ross. Akad. Nauk Sib. Otd. Inst. Mat., Novosibirsk, 1997,
  pp.~91--109, 190.

\bibitem[Ger98]{Gerasimov:fixedpoint}
V.N. Gerasimov, \emph{Fixed-point-free actions on cubings}, Siberian Adv. Math.
  \textbf{8} (1998), no.~3, 36--58.

\bibitem[GH07]{GuetnerHigson:catO}
E.~Guentner and N.~Higson, \emph{{Weak amenability of CAT(0)-cubical groups}},
  preprint, arXiv:math.OA/0702568, 2007.

\bibitem[Gui72]{Guichardet:bull}
A.~Guichardet, \emph{{Sur la cohomologie des groupes topologiques. II.}}, Bull.
  Sci. Math. \textbf{96} (1972), 305–332.

\bibitem[Gui77]{Guichardet}
\bysame, \emph{{\'Etude de la $1$-cohomologie et de la topologie du dual pour
  les groupes de Lie \`a radical ab\'elien}}, Math. Ann. \textbf{228} (1977),
  no.~1, 215--232.

\bibitem[Gui05]{Guirardel:trees}
V.~Guirardel, \emph{{Actions of finitely generated groups on R-trees}},
  preprint, arXiv: math.GR/0607295, 2005.


\bibitem[HP98]{HaPa}
F.~Haglund and F.~Paulin, \emph{{Simplicit\'e de groupes d'automorphismes
  d'espaces \`a courbure n\'egative}}, The Epstein Birthday Schrift (C.~Rourke,
  I.~Rivin, and C.~Series, eds.), Geometry and Topology Monographs, vol.~1,
  International Press, 1998, pp.~181--248.

\bibitem[HWK99]{Hagauer}
J.~Hagauer, W.Imrich, and S.~Klav\v{z}ar, \emph{Recognizing median graphs in
  subquadratic time}, Theoretical Computer Science \textbf{215} (1999),
  123--136.


\bibitem[Hou74]{Houghton:ends}
C.~H. Houghton, \emph{Ends of locally compact groups and their coset spaces},
  J. Austral. Math. Soc. \textbf{17} (1974), 274--284, Collection of articles
  dedicated to the memory of Hanna Neumann, VII.

\bibitem[Isb80]{Isbell}
J.~R. Isbell, \emph{Median {A}lgebra}, Trans. Amer. Math. Soc. (1980), no.~260,
  319--362.

\bibitem[JR06]{JohnsonRandrian}
W.~B. Johnson and N.~L. Randrianarivony, \emph{{$l\sb p\ (p>2)$} does
  not coarsely embed into a {H}ilbert space}, Proc. Amer. Math. Soc.
  \textbf{134} (2006), no.~4, 1045--1050.

\bibitem[KPR84]{KaltonPeck}
N.J. Kalton, N.T. Peck, and J.W. Roberts, \emph{{An $F$-space sampler}},
  {London Math. Soc. Lecture Note Series}, vol.~89, Cambridge University Press,
  Cambridge, 1984, xii+240 pp.

\bibitem[Laf]{Lafforgue:Trenforce}
V.~Lafforgue, \emph{Un renforcement de la propri\'et\'e {(T)}}, preprint, 2006.

\bibitem[MMR]{MMR}
F.R. McMorris, H.M. Mulder, and F.S. Roberts, \emph{{The median procedure on
  median graphs}}, Report 9413/B, Econometric Institute, Erasmus University
  Rotterdam.

\bibitem[Mul80]{Mulder:book}
H.~M. Mulder, \emph{The interval function of a graph}, Mathematical Centre
  Tracts, vol. 132, Mathematisch Centrum, Amsterdam, 1980.

\bibitem[Nac65]{Nachbin:Haar}
L.~Nachbin, \emph{{The Haar Integral}}, van Nostrand, Princeton, New Jersey,
  1965.

\bibitem[NR97]{NibloReeves:groupsactingCATO}
G.~A. Niblo and L.~D. Reeves, \emph{{Groups acting on CAT(0) cube complexes}},
  Geometry and Topology \textbf{1} (1997), 7 pp.

\bibitem[NR98]{NibloRoller:cubesandT}
G.~A. Niblo and M.~Roller, \emph{{Groups acting on cubes and Kazhdan's property
  (T)}}, Proc. Amer. Math. Soc. \textbf{126} (1998), no.~3, 693--699.

\bibitem[NSSS05]{NibloSageevScottSw}
G.~A. Niblo, M.~Sageev, P.~Scott, and G.A. Swarup, \emph{Minimal cubings}, Int.
  J. Algebra Comput. \textbf{15} (2005), no.~2, 343--366.

\bibitem[Nic]{NicaMaster}
B.~Nica, \emph{Group actions on median spaces}, Master Thesis.

\bibitem[Nic04]{Nica:cubulating}
\bysame, \emph{Cubulating spaces with walls}, Alg. \& Geom. Topology \textbf{4}
  (2004), 297--309.

\bibitem[Nie78]{Nieminen}
J.~Nieminen, \emph{The ideal structure of simple ternary algebras}, Colloq.
  Math. \textbf{40} (1978), 23--29.

\bibitem[Now06]{PNowak:lp}
P.W. Nowak, \emph{{Group actions on Banach spaces and a geometric
  characterization of a-T-menability}}, {Topology and its Applications}
  \textbf{153} (2006), 3409--3412.

  \bibitem[Now09]{PNowak:arxiv}
P.W. Nowak, \emph{{Group actions on Banach spaces and a geometric
  characterization of a-T-menability}}, preprint, {\texttt{http://arxiv.org/abs/math/0404402v3}}.

\bibitem[OP07]{OzawaPopa}
N.~Ozawa and S.~Popa, \emph{{On a class of II$_1$ factors with at most one
  Cartan subalgebra}}, preprint \textsc{arXiv:math.OA/07063623}, 2007.

\bibitem[Pan95]{Pansu:cohomologie}
P.~Pansu, \emph{{Cohomologie $L^p$: invariance sous quasi-isom\'etrie}},
  preprint, {\texttt{http://www.math.u-psud.fr/\%7Epansu/liste-prepub.html}},
  1995.

\bibitem[Rob98]{Robertson:Crofton}
G.~Robertson, \emph{{Crofton formulae and geodesic distance in hyperbolic
  spaces}}, J. Lie Theory \textbf{8} (1998), 163--172.


\bibitem[RS98]{RobertsonSteger:negdefker}
G.~Robertson and T.~Steger, \emph{{Negative definite kernels and a dynamical
  characterization of property (T) for countable groups}}, Ergod. Th. \& Dynam.
  Sys. \textbf{618} (1998), 247--253.


\bibitem[Rol98]{Roller:median}
M.~Roller, \emph{{Poc Sets, Median Algebras and Group Actions}}, preprint,
  University of Southampton,
  \texttt{http://www.maths.soton.ac.uk/pure/preprints.phtml}, 1998.

\bibitem[Sag95]{Sageev:ends}
M.~Sageev, \emph{{Ends of group pairs and non-positively curved cube
  complexes}}, Proc. London Maths. Soc. \textbf{71} (1995), no.~3, 585--617.


\bibitem[Sch85]{Scharlau85}
R.~Scharlau, \emph{{A characterization of Tits buildings by metrical
  properties}}, J. London Math. Soc. \textbf{32} (1985), no.~2, 317--327.


\bibitem[Sch38]{Schoenberg:negdefker}
I.J. Schoenberg, \emph{Metric spaces and positive definite functions}, Trans.
  Amer. Math. Soc. \textbf{44} (1938), 522--536.

\bibitem[Sel97]{Sela:acces}
Z.~Sela, \emph{{Acylindrical accessibility for groups}}, Invent. Math.
  \textbf{129} (1997), no.~3, 527--565.

\bibitem[Sho54a]{Sholander1}
M.~Sholander, \emph{Medians and betweenness}, Proc. Amer. Math. Soc. \textbf{5}
  (1954), 801--807.

\bibitem[Sho54b]{Sholander2}
\bysame, \emph{Medians, lattices and trees}, Proc. Amer. Math. Soc. \textbf{5}
  (1954), 808--812.


\bibitem[vdV84]{VandeVel}
M.~L.~J. van~de Vel, \emph{Dimension of binary convex structures}, Proc. London Math.
  Soc. \textbf{48} (1984), no.~3, 24--54.

\bibitem[vdV93]{VandeVel:book}
M.~L.~J. van~de Vel, \emph{Theory of convex structures}, North-Holland
  Mathematical Library, vol.~50, North-Holland Publishing Co., Amsterdam, 1993.

\bibitem[Ver93]{Verheul:book}
E.~R. Verheul, \emph{Multimedians in metric and normed spaces}, CWI Tract,
  vol.~91, Stichting Mathematisch Centrum voor Wiskunde en Informatica,
  Amsterdam, 1993.

\bibitem[Wat82]{Watatani}
Y.~Watatani, \emph{{Property (T) of Kazhdan implies property (FA) of Serre}},
  Math. Japonica \textbf{27} (1982), 97–--103.


\bibitem[WW75]{WellsWilliams:Embeddings}
J.H. Wells and L.R. Williams, \emph{Embeddings and extensions in analysis},
  Springer-Verlag, New York-Heidelberg, 1975.

\bibitem[Wil08]{Wildstrom}
J.~Wildstrom, \emph{Cost thresholds for dynamic resource location}, Discrete
  Applied Mathematics \textbf{156} (2008), 1846 –-- 1855.


\bibitem[Yu05]{Yu:hyplp}
G.~Yu, \emph{Hyperbolic groups admit proper isometric actions on
  $\ell^p$-spaces}, Geom. Funct. Anal. \textbf{15} (2005), no.~5, 1144--1151.

\end{thebibliography}
%%%%%%%%%%%%%%%%%%%%%%%%%%%%%%%%

%\end{document}

\newcommand{\etalchar}[1]{$^{#1}$}
\providecommand{\bysame}{\leavevmode\hbox to3em{\hrulefill}\thinspace}
\providecommand{\MR}{\relax\ifhmode\unskip\space\fi MR }
\providecommand{\MRhref}[2]{%
  \href{http://www.ams.org/mathscinet-getitem?mr=#1}{#2}
}
\providecommand{\href}[2]{#2}

\end{document}